\documentclass{article}
\usepackage{graphicx}
\usepackage[most]{tcolorbox}
\usepackage{amsmath}
\usepackage{amsthm}
\usepackage{amsfonts}
\usepackage{amssymb}
\usepackage{titlesec}
\usepackage{float}
\usepackage{bm}
\usepackage[mathcal]{euscript}
\usepackage{mathrsfs}
\usepackage[utf8]{inputenc}
\usepackage[font={small,it}]{subcaption}
\usepackage{hyperref}
\usepackage{geometry}
\usepackage[backend=biber,
            style=numeric,
            giveninits=true
            ]{biblatex}
\newtheorem{lemma}{Lemma}[section]
\newtheorem{remark}{Remark}[section]
\newtheorem{theorem}{Theorem}[section]

\numberwithin{equation}{section}

\newcommand{\ltwopm}{L^2(\Omega^{\pm})} 
\newcommand{\honepm}{H^1(\Omega^{\pm})}
\newcommand{\ltwo}{L^2(\Omega)}  
\newcommand{\bltwopm}{\bm{L}^2(\Omega^{\pm})}
\newcommand{\bhonepm}{\bm{H}^1(\Omega^{\pm})}
\newcommand{\bltwoc}{\bm{L}^2(\mathcal{C})}

\newcommand{\pc}{\pi}
\newcommand{\pch}{\pi_n^h}
\newcommand{\qc}{\theta}
\newcommand{\qch}{\theta^{h}}
\newcommand{\dotpc}{\dot{\pi}}
\newcommand{\ddotpc}{\ddot{\pi}}
\newcommand{\pcnh}{\pi_N^h}
\newcommand{\pcn}{\pi_n}
\newcommand{\pcnp}{\pi_{n-1}^h}
\newcommand{\pcnull}{\pi_{0}}
\newcommand{\frpcnh}{\delta_t \pi_{n}^h}
\newcommand{\enpc}{e_{n}^{\pi}}
\newcommand{\inpc}{a_{n}^{\pi}}
\newcommand{\enppc}{e_{n-1}^{\pi}}

\newcommand{\sfrac}{\psi^{c}}

\newcommand{\qhtr}{q^{h*}}
\newcommand{\pntr}{p_n^*}
\newcommand{\pnhtr}{p_n^{h*}}

\newcommand{\bu}{\bm{u}}
\newcommand{\dotp}{\dot{p}}
\newcommand{\dotu}{\dot{\bm{u}}}

\newcommand{\bv}{\bm{v}}
\newcommand{\bk}{\bm{K}}
\newcommand{\mc}{\mathcal{C}}
\newcommand{\telta}{\overline{\delta}}

\newcommand{\dunh}{\delta_{t}\bm{u}_n^h}
\newcommand{\dpnh}{\delta_{t}p_n^h}
\newcommand{\tnabla}{\nabla_{\tau}}

\newcommand{\pfs}{p_{n,k+1}^{h}}
\newcommand{\ppfs}{p_{n,k}^{h}}
\newcommand{\ufs}{\bm{u}_{n,k+1}^{h}}
\newcommand{\upfs}{\bm{u}_{n,k}^{h}}
\newcommand{\pcfs}{\pi_{n,k+1}^{h}}
\newcommand{\pfstr}{p_{n,k+1}^{h*}}
\newcommand{\pcpfs}{\pi_{n,k}^{h}}
\newcommand{\wfs}{w_{n,k+1}^{h}}
\newcommand{\wpfs}{w_{n,k}^{h}}
\newcommand{\tilu}{\widetilde{\bm{u}}_n^h}
\newcommand{\tilp}{\widetilde{p}_n^h}

\newcommand{\frun}{\delta_t \bm{u}_n}
\newcommand{\frunh}{\delta_t \bm{u}_n^h}
\newcommand{\frpnh}{\delta_t p_n^h}

\newcommand{\franu}{\delta_t a_n^{\bm{u}}}
\newcommand{\frenu}{\delta_t e_n^{\bm{u}}}

\newcommand{\enu}{e_n^{\bm{u}}}
\newcommand{\inu}{a^{\bm{u}}_n}
\newcommand{\inp}{a^p_n}
\newcommand{\enp}{e_n^p}
\newcommand{\enptr}{e_n^{p*}}
\newcommand{\inptr}{a_n^{p*}}

\newtcolorbox{strongproblem}[1][]{
    colback=white,
    colframe=black,
    boxrule=0.5pt,
    sharp corners,
}

\newtcolorbox{problem}[1][]{
    colback=white,
    colframe=black,
    boxrule=0.5pt,
    sharp corners,
    title=Problem (Q),
}

\newtcolorbox{discreteproblem}[1][]{
    colback=white,
    colframe=black,
    boxrule=0.5pt,
    sharp corners,
    title=Problem (Q-D),
}

\newtcolorbox{simpleproblem}[1][]{
    colback=white,
    colframe=black,
    boxrule=0.5pt,
    sharp corners,
    title=Problem (Q-D-S),
}

\title{Numerical analysis of coupled flow, deformation and frictional contact mechanics in a fractured porous medium}
\author{Marius Nevland%
\thanks{Center for Modeling of Coupled Subsurface Dynamics, Department of Mathematics, University of Bergen, Allégaten 41, Bergen, 5020, Norway.}\
\thanks{Corresponding author: marius.nevland@uib.no}\ ,
Kundan Kumar\footnotemark[1]\ ,
Inga Berre\footnotemark[1]\ ,
Jakub Wiktor Both\footnotemark[1]\ ,
Eirik Keilegavlen\footnotemark[1]
}
\date{}

\begin{document}

\setlength{\parindent}{0em}
\setlength{\parskip}{1em}
\maketitle

\begin{abstract}
This paper studies the numerical approximation of a mathematical model of coupled flow and viscoelastic deformation in a fractured porous medium including fracture contact mechanics. We consider a mixed-dimensional model, where the fracture is modeled as a co-dimension one object within the porous medium. This model couples Darcy flow in the matrix and fracture, linear viscoelasticity of the matrix and frictional contact mechanics of the fracture surfaces; the contact conditions include normal compliance and a Coulomb friction law. We propose and analyze a fully discrete numerical scheme for this problem, using conformal finite elements in space and the implicit Euler method in time. Existence and uniqueness of the discrete solution is established, and stability and a priori error estimates are derived. A numerical experiment is performed in which numerical error estimates are computed and compared to the theoretical results.
\end{abstract}
\renewcommand{\thefootnote}{}


\section{Introduction}

This paper is concerned with a mathematical model of coupled flow, deformation and frictional contact mechanics in a fractured porous medium. Such models find usage in several subsurface applications, such as CO$_2$-storage and geothermal energy production. The complexity of the coupled physical processes involved, together with the discontinuous nature of fractures, is reflected in the mathematical models being highly challenging and nonlinear.

Due to the typically high aspect ratio of fractures, it is common to model fractures as co-dimension one objects within the porous domain, leading to so-called mixed-dimensional models \cite{boon2021functional,martin}. Mixed-dimensional models of coupled flow, deformation and contact mechanics typically contain three components. First, the matrix (i.e. the porous medium surrounding the fractures) is modeled as linearly poroelastic, with fluid flow governed by Darcy's law; here, the textbook of Coussy \cite{coussy} is a standard reference. Second, Darcy flow is considered also in and across the fracture, using an averaged model to account for the collapsed dimension of the fracture \cite{martin}. The third component is the fracture contact mechanics, imposing constraints on the movement of the two sides of the fracture in the normal and/or tangential direction. The normal constraint is often formulated as a Karush-Kuhn-Tucker type inequality \cite{kikuchi_oden}, in which impenetrability between the two sides of the fracture is enforced (also called the Signorini condition). An alternative model is the normal compliance model, in which there is a constitutive relationship between the normal traction and penetration \cite{martins_oden}. The typical constraint in the tangential direction is the Coulomb friction model, in which the friction bound is proportional to the normal traction \cite{kikuchi_oden}. There have been several recent works on mathematical models containing some or all of the above components, focusing either on discretization approaches \cite{beaude2023mixed, berge, droniou_masson, garipov2016, GARIPOV2019104075, Stefansson_2021} or practical numerical algorithms for solving the discretized systems \cite{franceschini2022scalable,FRANCESCHINI2020113161, nevland2026augmented,zabegaev_efficient}.

Performing a rigorous numerical analysis of such models of coupled flow, deformation and fracture contact mechanics is a challenging task, due to the complexity of the models. Numerical analyses of a poroelastic model in a fractured porous medium, but excluding fracture contact mechanics, have been conducted in \cite{poroelastic_notmixed,Girault2019}, using conformal and mixed finite elements, respectively. Extending such analyses to include fracture contact mechanics is nontrivial, due to the nonlinear nature of the contact conditions. Nevertheless, there are some existing results in this direction. In an article of Bonaldi et al. \cite{bonaldi_multiphase}, a numerical scheme for a two-phase, mixed-dimensional poromechanical model, together with the Signorini condition and Coulomb friction law, was considered, and the existence and stability of a solution to the scheme was proven. This particular combination of contact conditions is known to cause significant mathematical difficulties \cite{chouly2023finite}, and as a result, the convergence to the continuous solution as the mesh size and time steps tend to zero was not proven. Another work of Bonaldi et al. \cite{bonaldi2024numerical} considered a numerical scheme for mixed-dimensional poromechanics with a Signorini condition on the fracture, but without friction, and here the convergence of the discrete solution to the continuous solution in the limit was proven, in addition to existence and stability estimates. Neither of the works included a priori error estimates between the continuous and discrete solutions.

In our previous work \cite{nevland2026numerical}, we analyzed a space-time discretization of the quasi-static Biot equations for a poroviscoelastic medium subject to frictional contact with a rigid obstacle; our analysis included proving the existence and uniqueness of the discrete solution, stability and a priori error estimates. In the present paper, we extend the analysis of \cite{nevland2026numerical} to a mixed-dimensional model of poroviscoelasticity in a fractured medium, including a linearized fracture flow model as well as fracture contact mechanics. Our fracture flow model is based on the work of Martin et al. \cite{martin}, in which a dimensionally reduced mass balance equation is derived by averaging across the width of the fracture. In particular, this model accounts for discontinuous pressures across matrix-fracture interfaces. Our contact mechanics model is the same as in our previous paper \cite{nevland2026numerical}, but extended to take into account the coupling to fracture pressure, and consists of the normal compliance and Coulomb friction laws.

The numerical scheme employed in the present paper consists of conformal finite elements in space and the implicit Euler method in time. Existence, uniqueness and stability of the discrete solution is established, and a priori error estimates and corresponding convergence rates are proven. To our knowledge, this is the first time that combined space-time error estimates are derived for a mixed-dimensional model of coupled flow, deformation and contact mechanics. The existence is proven by a method first established in our previous paper \cite{nevland2026numerical}, extended here to the mixed-dimensional setting. This method consists of first considering the simpler Tresca friction model, and proving the existence to this simpler problem by decoupling the flow and mechanics by the fixed-stress splitting scheme. A fixed-point argument is then used to arrive at the solution to the problem with Coulomb friction. The stability and error estimates are derived by combining techniques from the aforementioned works on mixed-dimensional poroelasticity without contact mechanics \cite{poroelastic_notmixed,Girault2019} and works on purely mechanical contact problems \cite{hemivariational_ineq,han2002quasistatic}. The theoretical analysis is supplemented by a numerical experiment, in which numerical convergence orders are computed. It is observed that the reduced regularity of the solution around fracture tips adversely affects the convergence rate. 

The paper is structured as follows. We start by introducing the mixed-dimensional mathematical model in Section \ref{sec:model}, and the variational formulation of the model is derived in Section \ref{sec:weak_form}. In Section \ref{sec:discretization}, our fully discrete numerical scheme is proposed and analyzed. The numerical experiment is presented in Section \ref{sec:num_exp}. Finally, concluding remarks are provided in Section \ref{sec:conclusion}.

\section{Model formulation} \label{sec:model}

\subsection{Notation and geometrical setting}
\begin{figure}
[H]\centering
    {\scalebox{0.25}{\includegraphics{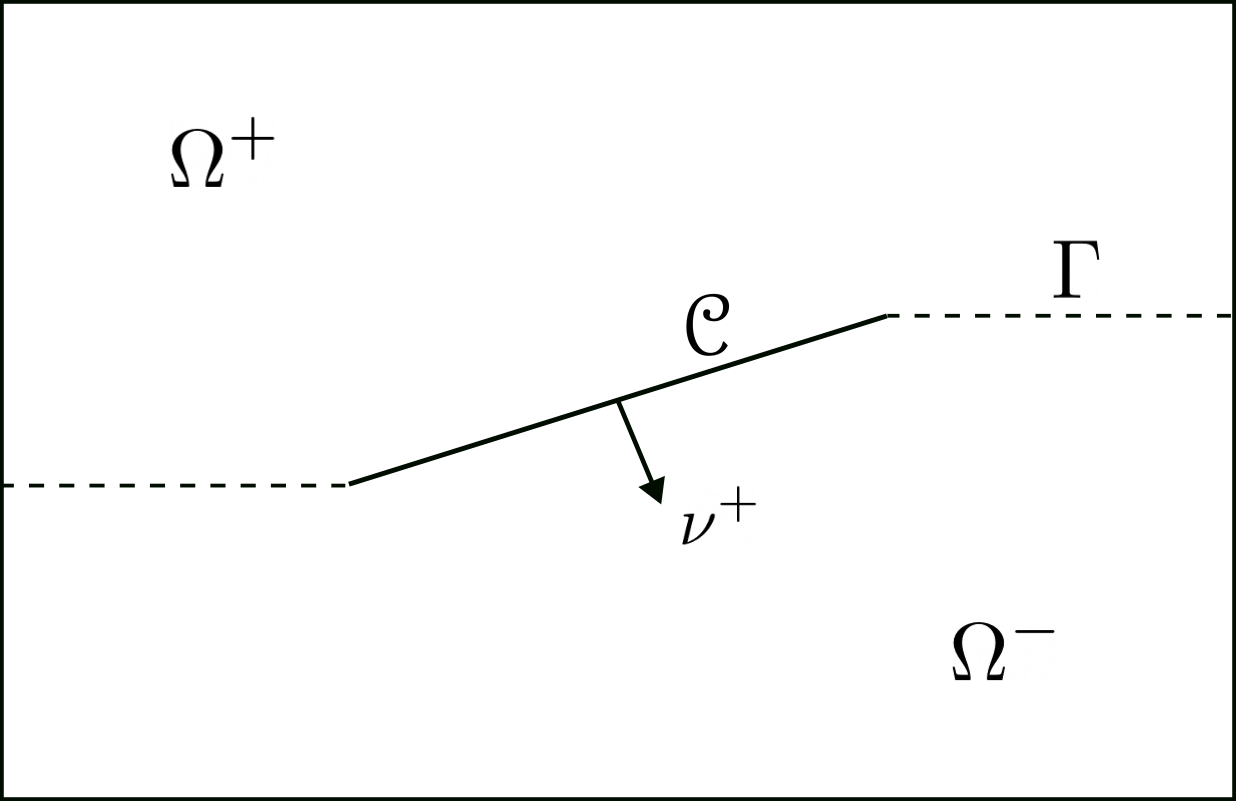}}}
    \caption{Schematic of the domain, fracture and interface.}
    \label{fig:domain}
\end{figure}

Let the reservoir $\Omega \subset \mathbb{R}^d$ be a open, bounded domain, $d=2$ or 3, with piecewise smooth, Lipschitz boundary $\partial \Omega$ and exterior normal $\bm{\nu}$. Denote by $\overline{\Omega}$ the closure of the domain. The fracture $\mathcal{C} \subset \Omega$ is taken to be a line segment when $d=2$ and a plane segment when $d=3$. The domain $\Omega$ is separated into two non-overlapping, open subdomains $\Omega^+$ and $\Omega^-$ with a Lipschitz interface $\Gamma$ containing $\mathcal{C}$. Thus, we have $\Omega=\Omega^+ \cup \Omega^- \cup \Gamma$. See Figure \ref{fig:domain} below for a visualization of the domain.
We shall also make frequent use of the notation: 
\begin{equation}
    \Omega^{\pm} := \Omega^+ \cup \Omega^-.
\end{equation}
The superscript $*$ will be used to denote either $+$ or $-$. We denote by $\bm{\nu}^*$ the outward normal vector of $\Omega^*$. As the two subdomains $\Omega^+$ and $\Omega^-$ share an interface, we have the relation $\bm{\nu}^-=-\bm{\nu}^+$ on $\Gamma$. 

For a function $f$ in $\Omega^{\pm}$ with a well-defined trace, we let $f^*$ denote its trace from $\Omega^*$ on $\Gamma$. The jump of $f$ over $\Gamma$ is defined by
\begin{equation} \label{eq:jump_operator}
    [f]=f^+-f^-.
\end{equation}
For vector-valued functions $\bm{v}$ with a well-defined trace, we define the normal and tangential components of its jump over $\Gamma$ as:
\begin{equation} \label{eq:jump_norm_tang}
    [\bm{v}]_{\nu}=[\bm{v}] \cdot \bm{\nu}^+ \ , \ [\bm{v}]_{\tau}=[\bm{v}]-[\bm{v}]_{\nu}\bm{\nu^+}.
\end{equation}
The time derivative of a function $f$ will mostly be denoted by a dot, $\dot{f}$, and its double derivative by a double dot, $\ddot{f}$. Occasionally, we will instead use the notation $\partial_t$ to denote a time derivative.

\subsection{Matrix poroviscoelasticity}

In the matrix subdomain, we model coupled flow and deformation by the quasi-static Biot equations for a viscoelastic, homogeneous, isotropic, porous solid saturated with a compressible viscous fluid. The equations are formulated in terms of the unknown variables pressure $p$ and displacement $\bu$. The mathematical model in the matrix is formulated by a domain decomposition approach, where the equations are imposed on $\Omega^{\pm}$, while separate conditions ensure the continuity of the solution across the part of the interface not containing the fracture, $\Gamma \setminus \mc$.

Balance of momentum of the solid reads
\begin{equation} \label{eq:mom_bal}
    -\nabla \cdot \bm{\sigma}^{\text{por}}=\bm{f} \ \ \ \text{in} \ \ \Omega^{\pm} \times (0,T),
\end{equation}
where $\bm{\sigma}^{\text{por}}$ denotes the total poroviscoelastic stress tensor and $\bm{f}$ denotes body forces. The constitutive relation for $\bm{\sigma}^{\text{por}}$ is
\begin{equation} \label{eq:stress_total}
    \bm{\sigma}^{\text{por}}(\bm{u},p)=\bm{\sigma}(\bm{u})-\alpha p \bm{I},
\end{equation}
where $\bm{I}$ is the identity tensor,  $\alpha>0$ is the dimensionless Biot coefficient, and $\bm{\sigma}$ is the viscoelastic stress tensor, following the Kelvin-Voigt model \cite{Sofonea_Matei_2012}:
\begin{equation} \label{stress_mechanics}
    \bm{\sigma}(\bm{u})=2G\bm{\varepsilon}(\bm{u})+\lambda(\nabla \cdot \bm{u})\bm{I}+\gamma\bm{\varepsilon}(\dotu)+\gamma(\nabla \cdot \dotu)\bm{I},
\end{equation}
where $\lambda$ and $G$ are the Lamé constants , $\gamma>0$ is a viscoelastic parameter, and $\bm{\varepsilon}(\bm{u})$ is the linearized strain tensor:
\begin{equation} \label{eq:strain_tensor}
    \bm{\varepsilon}(\bm{u})=\frac{1}{2}\left(\nabla\bm{u}+(\nabla\bm{u})^T
\right).
\end{equation}
The volumetric fluid flux follows Darcy's law (for simplicity, we neglect gravitational effects):
\begin{equation} \label{darcy}
    \bm{v}_D=-\frac{\bm{K}}{\mu}\nabla p,
\end{equation}
where $\bm{K}$ is the permeability tensor, and $\mu>0$ is the constant fluid viscosity. The permeability tensor is assumed to be time-independent, symmetric, bounded, and uniformly positive definite in space; that is, there exists $\lambda_{\text{min}}>0, \lambda_{\text{max}}>0$ such that for all $\bm{x}\in \Omega$:
\begin{equation} \label{eq:spd_perm}
    \lambda_{\text{min}}\lvert \bm{\xi}\rvert^2 \leq \bm{\xi}^T\bm{K}(\bm{x})\bm{\xi} \leq \lambda_{\text{max}}\lvert \bm{\xi}\rvert^2.
\end{equation}
The fluid is assumed to be slightly compressible, and we consider the simple case when the fluid density $\rho$ is a linear function of the pressure:
\begin{equation} \label{eq:density_pressure}
    \rho=\rho_{0}(1+c_f(p-p_0)),
\end{equation}
where $c_f$ is the fluid compressibility, and subscripts 0 denote reference values of the respective variables. The porosity $\phi$ is given by
\begin{equation} \label{eq:fluid_content}
    \phi=\phi_0+\alpha\nabla \cdot \bm{u}+\frac{1}{M}p,
\end{equation}
where $\phi_0$ is the initial reference porosity and $M$ is a Biot constant.
The fluid mass balance equation reads
\begin{equation} \label{eq:mass_balance}
\partial_t(\rho\phi)+\nabla \cdot (\rho \bm{v}_D)=\psi,
\end{equation}
with $\psi$ denoting a mass source or sink term. This equation is subsequently linearized, by neglecting small quantities, following \cite{poroelastic_notmixed}. Specifically, by expanding the time derivative in \eqref{eq:mass_balance} and collecting factors of $\dotp$, we obtain, after substituting the constitutive relations \eqref{darcy},\eqref{eq:density_pressure} and \eqref{eq:fluid_content}:
\begin{equation} \label{eq:biot_intermediate}
    \begin{aligned}
    &\rho_0\left(\frac{1}{M}(1+c_f(p-p_0))+c_f(\phi_0+\alpha\nabla \cdot \bu + \frac{1}{M}p)\right)\partial_tp\\
    &+\rho_0(1+c_f(p-p_0))\alpha\nabla \cdot \partial_t\bu+\nabla \cdot (\rho_0(1+c_f(p-p_0)\bv_D)=\psi.
    \end{aligned}
\end{equation}
We now assume that $c_f(p-p_0), p/M$ and $\nabla \cdot \bu$ are all small quantities. Hence, we approximate the density $\rho_0(1+c_f(p-p_0))$ by $\rho_0$ everywhere in \eqref{eq:biot_intermediate}, as well as approximating $c_f(\phi_0+\alpha\nabla \cdot \bu + p/M)$ by $c_f\phi_0$. For simplicity, we also define $c_0=1/M+c_f\phi_0$. With these approximations, the resulting linearized mass balance equation becomes, after dividing by $\rho_0$:
\begin{equation} \label{eq:biot_mass_balance}
  c_0 \dotp+\alpha\nabla \cdot \dotu-\frac{1}{\mu}\nabla \cdot (\bk\nabla p)=\psi \ \ \ \text{in} \ \ \Omega^{\pm} \times (0,T),
\end{equation}
where by abuse of notation, we denote the rescaled source term again by $\psi$.

\subsection{Fracture mass balance}

The width of the fracture is represented by a function $w$ on $\mathcal{C}$, defined as the sum of the reference fracture width in the undeformed configuration, denoted by $g$, and the jump of $\bm{u} \cdot \bm{\nu}^-$ on $\mathcal{C}$:
\begin{equation} \label{eq:fracture_width}
    w=-[\bm{u}]_{\nu}+g.
\end{equation}
We assume the reference fracture width $g$ to be bounded from both above and below, that is, for all $\bm{x}\in \mc$:
\begin{equation} \label{eq:g_bound}
    0<g_{\text{min}} \leq g(\bm{x})\leq g_{\text{max}}.
\end{equation}
In particular, this implies that the fracture is mechanically open in the reference configuration. Our fracture fluid flow model is similar to the one considered in \cite{stefansson2023flexible}, which is in turn inspired by the work of \cite{martin}, where a dimensionally reduced fracture flow model was derived by averaging across the width of the fracture. This averaging procedure results in the width $w$ appearing as a scaling factor in the terms. Denote by $\pc$ the fracture pressure, and denote by $\tnabla$ the gradient on $\mathcal{C}$, acting only along directions tangential to $\mc$. The volumetric fluid flux is modeled by Darcy's law with a cubic law for the permeability:
\begin{equation} \label{eq:darcy_fracture}
    \bm{v}_{\mathcal{C}}=-\frac{w^2}{12\mu}\tnabla \pc.
\end{equation}
The dimensionally reduced mass balance equation in the fracture reads:
\begin{equation} \label{eq:mass_balance_fracture_full}
    \partial_t(w\rho_{\mc})+\tnabla \cdot (w\rho_{\mc}\bm{v}_{\mathcal{C}})-\sum_{*\in \pm} \rho^{\text{up}} \bm{v}_D^*\cdot \bm{\nu}^*=w\sfrac,
\end{equation}
where $\sfrac$ is a known mass source or sink term, $\rho^{\text{up}}$ is an upstream-weighted ``interface density" (defined below), and $\rho_{\mc}$ is the fluid density in the fracture, depending on the fracture pressure by the same constitutive law as before:
\begin{equation} \label{eq:density_fracture}
    \rho_{\mc}=\rho_0(1+c_f(\pc-\pc_0)).
\end{equation}
The third term on the left-hand side of \eqref{eq:mass_balance_fracture_full} is a coupling term that accounts for fluid mass flux between the matrix and fracture, given as the sum of the normal traces of the matrix mass flux from each side of the fracture. Following \cite{martin}, the normal trace of the Darcy flux, $\bm{v}_D^* \cdot \bm{\nu}^*$, is required to satisfy the following averaged Darcy's law:
\begin{equation} \label{eq:normal_darcy}
    \bm{v}_D^*\cdot \bm{\nu}^*=-\frac{2K_{\nu}}{w\mu}(\pc-p^*)
\end{equation}
where $K_{\nu}$ is the normal permeability, which is also modeled by a cubic law:
\begin{equation} \label{eq:normal_perm}
    K_{\nu}=\frac{w^2}{12}.
\end{equation}
Finally, we define the upstream-weighted density $\rho^{\text{up}}$, which equals either the trace of the matrix density or the fracture density, depending on the sign of $\bm{v}_D^* \cdot \bm{\nu}^*$:
\begin{equation}
    \rho^{\text{up}}=\begin{cases}
        \rho^* \ \ \text{if} \ \bm{v}_D^* \cdot \bm{\nu}^* >0,\\
        \rho_{\mc} \ \ \text{if} \  \bm{v}_D^* \cdot \bm{\nu}^* \leq 0.
    \end{cases}
\end{equation}
Some simplifications to equation \eqref{eq:mass_balance_fracture_full} will be made as follows: We expand the time derivative term, and use \eqref{eq:fracture_width} and the fact that $g$ is constant in time to get
\begin{equation}  \partial_t(w\rho_{\mc})=\rho_{\mc}\partial_tw+w\partial_t\rho_{\mc}=\rho_{\mc}\partial_t(-[\bm{u}]_{\nu})+w\partial_t \rho_{\mc}.
\end{equation}
Approximating $\rho_{\mc}$ by $\rho_0$ in the first term, under the assumption that $c_f(\pc-\pc_0)$ is small, and differentiating $\rho_{\mc}$ in the second term according to the expression \eqref{eq:density_fracture}, we get
\begin{equation}
    \partial_t(w\rho_{\mc})\approx -\rho_0\big([\dot{\bm{u}}]_{\nu}-w c_f \dotpc\big).
\end{equation}
We similarly approximate $\rho_{\mc}$ and $\rho^{\text{up}}$ by $\rho_0$ in the second and third terms of \eqref{eq:mass_balance_fracture_full}, respectively. The result, after substituting \eqref{eq:darcy_fracture}, \eqref{eq:normal_darcy} and \eqref{eq:normal_perm} and dividing by $\rho_0$, is the reduced fracture mass balance equation:
\begin{equation}
    -[\dot{\bm{u}}]_{\nu}+c_f w \dotpc-\frac{1}{\mu}\tnabla \cdot (\frac{w^3}{12}\tnabla \pc)+\sum_{*\in\pm}\frac{w}{6\mu}(\pc-p^*)=w\sfrac,
\end{equation}
where we reuse the notation $\sfrac$ to denote the rescaled source term. The equation is further linearized by approximating the width $w$ by the width in the undeformed configuration, $w\approx g$. Thus, our linearized fracture mass balance equation is:
\begin{equation} \label{eq:mass_balance_fracture}
   -[\dot{\bm{u}}]_{\nu}+c_f g \dotpc-\frac{1}{\mu}\tnabla \cdot (\frac{g^3}{12}\tnabla \pc)+\sum_{*\in\pm}\frac{g}{6\mu}(\pc-p^*)=g\sfrac \ \ \ \text{on} \ \ \mc \times (0,T). 
\end{equation}

\subsection{Interface, boundary and initial conditions} \label{sec:contact_mechanics}

As previously mentioned, continuity conditions are enforced on the part of the interface not including the fracture, $\Gamma\setminus\mc$. For the flow, we require the continuity of the pressure and normal flux across the interface:
\begin{equation} \label{eq:interface_cond_flow}
    [p]=0, \ [\bv_D]_{\nu}=0 \ \ \text{on} \ \ \Gamma\setminus\mc \times (0,T),
\end{equation}
and for the mechanics, we require the continuity of the displacement and the normal stresses across the interface:
\begin{equation} \label{eq:interface_cond_mech}
    [\bu]=0, \ [\bm{\sigma}^{\text{por}}]_{\nu}=\bm0 \ \ \text{on} \ \ \Gamma\setminus\mc \times (0,T).
\end{equation}
On the fracture $\mc$, we impose frictional contact conditions. We start by defining the effective contact tractions on the two fracture surfaces, denoted by $\bm{t}^+$ and $\bm{t}^-$. Noting that the fracture fluid pressure contributes to the effective contact tractions, they are defined as follows:
\begin{equation}
\begin{cases} \label{eq:contact_traction}
    \bm{t}^*=(\bm{\sigma}^{\text{por}})^* \bm{\nu}^*+\pc\bm{\nu}^*, \ \ * \in \pm, \\
    t_{\nu}=\bm{t}^+ \cdot \bm{\nu}^+, \\
    \bm{t}_{\tau}=\bm{t}^+-t_{\nu}\bm{\nu}^+.
\end{cases}
\end{equation}
By Newton's third law, we must have:
\begin{equation} \label{eq:interface_force_balance}
    \bm{t}^+ + \bm{t}^-=\bm{0} \ \ \ \text{on} \ \ \mc \times (0,T) .
\end{equation}
For the normal component of the contact traction, we impose a normal compliance law, similar to the one used by Martins and Oden \cite{martins_oden}. It introduces a constitutive relationship between the normal contact traction and the penetration, and we formulate it as follows:
\begin{equation} \label{eq:normal_compliance}
    t_{\nu} = -c_{\nu}([\bm{u}]_{\nu}-g)_+ \ \ \ \text{on} \ \ \mc \times (0,T),
\end{equation}
where $c_{\nu}>0$ is a constant, which can be interpreted as the normal stiffness of the fracture, and $(\cdot)_+$ denotes the positive cut of the function. 

The relative tangential movement of the fracture is constrained by the Coulomb friction law. This law states that during contact, slip occurs if the tangential contact traction reaches
a friction bound, while no relative movement happens if the bound is not reached.
Moreover, this friction bound is proportional to the magnitude of the normal contact traction, with the
proportionality constant being the friction coefficient, $\mu_{\text{fr}}$. Taking into account the normal compliance law \eqref{eq:normal_compliance} for the normal contact traction, the Coulomb friction law, imposed on $\mc \times (0,T)$, reads:
\begin{equation} \label{eq:friction_law}
\begin{aligned}
&\text{If} \ [\bu]_{\nu} \leq g \ \text{then} \ \ \lvert \bm{t}_{\tau}\rvert =0,\\
   & \text{If} \ [\bu]_{\nu} >g  \ \text{then} \ \begin{cases}
        \lvert \bm{t}_{\tau}\rvert \leq \mu_{\text{fr}} c_{\nu}([\bu]_{\nu}-g),\\
        \lvert \bm{t}_{\tau}\rvert < \mu_{\text{fr}} c_{\nu}([\bu]_{\nu}-g)  \implies [\dot{\bu}]_{\tau}=\bm{0}, \\
        \lvert \bm{t}_{\tau}\rvert = \mu_{\text{fr}} c_{\nu}([\bu]_{\nu}-g) \implies \exists \ \zeta \geq 0 \ , \ [\dot{\bu}]_{\tau}=-\zeta \bm{t}_{\tau}.
    \end{cases}
\end{aligned}
\end{equation}
Next, we consider conditions on the exterior boundary $\partial \Omega$. To this end, we introduce two different decompositions of $\partial \Omega$ into disjoint subsets, $\partial \Omega=\partial \Omega_p \cup \partial \Omega_f$ and $\partial \Omega=\partial \Omega_d \cup \partial \Omega_n$, corresponding to the flow and mechanical boundary conditions, respectively. Dirichlet conditions are prescribed on $\partial\Omega_d$ and $\partial\Omega_p$, and Neumann conditions on $\partial \Omega_f$ and $\partial \Omega_n$. The Dirichlet boundaries $\partial \Omega_d$ and $\partial \Omega_p$ are assumed to contain positive-measure subsets of both $\partial\Omega^+$ and $\partial\Omega^-$. To simplify the analysis, we consider homogeneous boundary conditions:
\begin{equation} \label{eq:bcs}
\begin{aligned}
    &p=0 \ \ \ \text{on} \ \ \partial \Omega_p\times (0,T), \ \  \bk\nabla p\cdot\bm{\nu}=0 \ \ \text{on} \ \ \partial \Omega_f\times (0,T),\\
    &\bu=\bm{0} \ \ \  \text{on} \ \ \partial \Omega_d\times (0,T), \ \ \bm{\sigma}^{\text{por}}\cdot\bm{\nu}=\bm{0} \ \ \text{on} \ \ \partial \Omega_n\times(0,T),
\end{aligned}
\end{equation}
Note, however, that in the numerical experiments of Section \ref{sec:num_exp}, a nonhomogeneous mechanical Neumann condition is used. Since a nonhomogeneous Neumann condition can be absorbed into the source terms $\bm{f}$ and/or $\psi$ in the variational formulation, there is no loss of generality in assuming homogeneous Neumann conditions for the analysis.

Additionally, we prescribe zero normal flux (homogeneous Neumann condition) on the boundary of $\mc$, i.e. the fracture tips:
\begin{equation} \label{eq:fracture_tips}
    \frac{g^3}{12}\tnabla \pc \cdot \bm{\nu}_{\mc}=0 \ \ \text{on} \ \ \partial\mc \times (0,T),
\end{equation}
where $\bm{\nu}_{\mc}$ denotes the outward normal on $\partial \mc$; when $d=2$, the condition reduces to just $\frac{g^3}{12}\tnabla\pc = 0$ at the two points constituting $\partial \mc$. Such a condition is supported by experimental evidence \cite{both2024}. Finally, initial values of the displacement and pressure are prescribed: 
\begin{equation} \label{eq:init_cond}
    \bu(0)=\bu_0, \  p(0)=p_0, \ \pc(0)=\pcnull
\end{equation}
Summarizing, our problem to be solved is the following:
\begin{strongproblem}
    Find $\bu: \Omega^{\pm} \times (0,T) \to \mathbb{R}^d$, $p: \Omega^{\pm} \times (0,T) \to \mathbb{R}$ and $\pc: \mc \times (0,T) \to \mathbb{R}$ satisfying equations \eqref{eq:mom_bal}, \eqref{eq:biot_mass_balance} and \eqref{eq:mass_balance_fracture}, subject to the interface continuity conditions \eqref{eq:interface_cond_flow}, \eqref{eq:interface_cond_mech}, the contact conditions \eqref{eq:interface_force_balance}, \eqref{eq:normal_compliance}, \eqref{eq:friction_law}, the boundary conditions \eqref{eq:bcs}, \eqref{eq:fracture_tips} and the initial conditions \eqref{eq:init_cond}.
\end{strongproblem}

\section{Variational formulation} \label{sec:weak_form}

In this section, the model of the previous section is put in a variational formulation. We first review the functional-analytic notation to be used in this paper.

\subsection{Notation and functional setting}

Recall our division of the domain into two non-overlapping subdomains $\Omega^+$ and $\Omega^-$, separated by an interface $\Gamma$ containing the fracture $\mc$. For quantities defined on $\Omega^{\pm}$, we will frequently take the sum of their $L^2$ inner products over $\Omega^+$ and $\Omega^-$. Denoting each of these individual inner products by $(\cdot,\cdot)_{\Omega^+}$ and $(\cdot,\cdot)_{\Omega^-}$, respectively, we denote their sum simply by parentheses without a subscript, i.e., 
\begin{equation} \label{eq:l2_ip_combined}
    (\cdot,\cdot)=(\cdot,\cdot)_{\Omega^+}+(\cdot,\cdot)_{\Omega^-}.
\end{equation}
A subscript $\mc$ is used to denote the $L^2$ inner product over the fracture, $(\cdot,\cdot)_{\mc}$. The $L^2$ inner product is used for functions with values in both $\mathbb{R}, \mathbb{R}^d$ and $\mathbb{S}^d$ (the space of second order symmetric tensors on $\mathbb{R}^d$); the above parenthesis notation is used in all three cases. The norm induced from the inner product \eqref{eq:l2_ip_combined} is denoted by the subscript $\ltwopm$:
\begin{equation}
    \lVert \cdot \lVert_{\ltwopm}=\left(\lVert \cdot \rVert_{L^2(\Omega^+)}+\lVert \cdot \rVert_{L^2(\Omega^-)}\right)^{1/2}.
\end{equation}
Since we are working with variables that may jump across the fracture, it is necessary to introduce broken Sobolev spaces. We denote by $H^s(\Omega^{\pm})$ the broken Sobolev space defined as follows:
\begin{equation} \label{eq:h1pm}
    H^s(\Omega^{\pm})=\{v\in L^2(\Omega); \ v\rvert_{\Omega^*}\in H^s(\Omega^*), \ *\in \pm \},
\end{equation}
for a real number $s \geq 0$. In the case when $s$ is an integer, $H^s(\Omega^*)$ denotes the usual Sobolev space of functions whose weak derivatives up to order $s$ lie in $L^2(\Omega^*)$. If $s$ is not an integer, $H^s(\Omega^*)$ denotes a fractional Sobolev space, which can be defined using interpolation theory of Banach spaces \cite{adams2003sobolev}. The broken Sobolev spaces are normed by
\begin{equation}
    \lVert v \rVert_{H^s(\Omega^{\pm})}=\left(\lVert v\rVert_{H^s(\Omega^+)}^2+\lVert v \rVert_{H^s(\Omega^-)}^2\right)^{1/2},
\end{equation}
where $\lVert \cdot \rVert_{H^s(\Omega^*)}$ denotes the usual Sobolev norm if $s$ is an integer, and otherwise, it denotes a fractional Sobolev norm; see \cite{adams2003sobolev} for details of the latter. Additionally, we denote by $\lvert \cdot \rvert_{H^s(\Omega^*)}$ the standard seminorm for these spaces. We also define the Sobolev spaces $H^s(\mc)$ for functions defined on the fracture, where it is understood that the derivatives act only along directions tangential to $\mc$. In cases where the output space is either $\mathbb{R}^d$ or $\mathbb{S}^d$, the corresponding Sobolev space is denoted by a boldface, for instance $\bm{H}^s(\Omega^{\pm})$, where the output space should be clear from the context.

The space for the displacements, denoted by $\bm{V}$, is defined as follows:
\begin{equation} \label{eq:func_space_displacement}
    \bm{V}=\{\bm{v} \in \bm{H}^1(\Omega^{\pm}) \ ; \ [\bm{v}]=\bm{0} \ \text{on} \ \Gamma\setminus\mc, \ \bm{v}\lVert_{ \partial\Omega_d}=\bm{0} \},
\end{equation}
equipped with the norm of $\bm{H}^1(\Omega^{\pm})$, namely:
\begin{equation}
  \lVert \bv \rVert_{\bhonepm}=\left(\sum_{i=1}^d \lVert v_i \rVert^2_{\honepm}\right)^{1/2}.  
\end{equation}
The space for the pressure in the matrix is the following:
\begin{equation} \label{eq:pressure_space_notmixed}
    Q=\{q\in \honepm \ ; \ [q]=0 \ \text{on} \ \Gamma\setminus\mc, \ q\lvert_{\partial\Omega_p}=0\},
\end{equation}
equipped with the norm of $\honepm$. Finally, the space for the pressure in the fracture is $H^1(\mc)$.

Next, we recall some important inequalities to be used in the analysis. Recall that the trace of a function $v \in H^1(\Omega^*)$ on $\mc$, is denoted by $v^* \in L^2(\mc)$. The usual trace inequality holds \cite{poroelastic_notmixed}:
\begin{equation} \label{eq:trace_inequality}
    \forall v \in H^1(\Omega^*), \ \ \ \lVert  v^* \rVert_{L^2(\mc)} \leq C^*\lVert v \rVert_{H^1(\Omega^*)}, \ *\in \pm,
\end{equation}
for some $C^*>0$ depending only on $\Omega^*$ and $\mc$. 

Since $\partial\Omega_p$ and $\partial\Omega_d$ contain positive-measure subsets of both $\partial\Omega^+$ and $\partial \Omega^-$, we may use Poincaré's inequality and Korn's first inequality separately on $\Omega^+$ and $\Omega^-$, and then combine them into inequalities on $\Omega^{\pm}$. The resulting Poincaré's inequality on the space $Q$ reads:
\begin{equation} \label{eq:poincare_inequality}
    \forall v \in Q, \ \ \lVert v \rVert_{L^2(\Omega)} \leq C\lvert v \rvert_{\honepm},
\end{equation}
for a constant $C$ independent of $v$; an equivalent inequality holds for $\bv \in \bm{V}$. Korn's first inequality on the space $\bm{V}$ reads:
\begin{equation} \label{eq:korn_inequality}
    \forall \bm{v} \in \bm{V}, \ \ \lvert \bm{v} \rvert_{\bhonepm} \leq C\lVert \bm{\varepsilon}(\bm{v})\rVert_{\bltwopm},
\end{equation}
for a constant $C$ independent of $\bm{v}$.

Time-dependent $L^p$-spaces (or Bochner spaces), for functions defined on a time interval $(a,b)$ with values in a space $X$, are defined as follows:
\begin{equation} \label{eq:bochner_space}
    L^p(a,b;X)=\left\{ f \ \text{measurable in} \ (a,b) \ ; \ \int_a^b\lVert f(t) \rVert_X^p dt < \infty\right\},
\end{equation}
equipped with the norm
\begin{equation} \label{eq:norm_bochner}
    \lVert f \rVert_{L^p(a,b;X)}=\left(\int_a^b \lVert f(t) \rVert_X^p dt\right)^{1/p},
\end{equation}
with the usual modification if $p=\infty$. We shall also need the spaces $H^1(a,b;X)$ and $H^2(a,b;X)$, defined as follows:
\begin{equation}
    H^1(a,b;X)= \{ f\in L^2(a,b;X) \ ; \ \dot{f} \in L^2(a,b;X)\}.
\end{equation}
\begin{equation}
    H^2(a,b;X)= \{ f\in H^1(a,b;X) \ ; \ \ddot{f} \in L^2(a,b;X)\}.
\end{equation}
equipped with the norms
\begin{equation} \label{eq:norm_h1_bochner}
    \lVert f \rVert_{H^1(a,b;X)}=\left(\int_a^b \lVert f(t) \rVert_X^2 + \lVert \dot{f}(t) \rVert^2_X dt\right)^{1/2},
\end{equation}
\begin{equation} \label{eq:norm_h2_bochner}
    \lVert f \rVert_{H^2(a,b;X)}=\left(\int_a^b \lVert f(t) \rVert_X^2 + \lVert \dot{f}(t) \rVert^2_X + \lVert \ddot{f}(t) \rVert^2_X dt\right)^{1/2}.
\end{equation}
Finally, we denote the dual spaces of $H^s$ by $H^{-s}$, and otherwise a prime is used to denote the dual space, e.g., $Q^{\prime}$ is the dual space of $Q$. Duality pairings are denoted by angular brackets, with the same subscript conventions as for the $L^2$ inner product. That is, $\langle \cdot,\cdot\rangle$ and $\langle \cdot,\cdot\rangle_{\mc}$ denote duality pairings between spaces defined over $\Omega^{\pm}$ and $\mc$, respectively.

\subsection{Variational formulation} \label{sec:variational_form}

In order to make the presentation of the variational formulation more compact, we first define some operators which are standard in the contact mechanics literature. The elastic and viscoelastic bilinear forms $A[\cdot,\cdot]: \bm{V} \times \bm{V} \to \mathbb{R}$ and $B[\cdot,\cdot]: \bm{V} \times \bm{V} \to \mathbb{R}$ are defined as follows
\begin{equation} \label{eq:bilinear_form_1}
    A[\bm{u},\bm{v}]=2G(\bm{\varepsilon}(\bm{u}),\bm{\varepsilon}(\bm{v}))+\lambda(\nabla \cdot \bm{u},\nabla \cdot \bm{v}),
\end{equation}
\begin{equation} \label{eq:bilinear_form_2}
B[\bm{u},\bm{v}]=\gamma(\bm{\varepsilon}(\bm{u}),\bm{\varepsilon}(\bm{v}))+\gamma(\nabla \cdot \bm{u},\nabla \cdot \bm{v}).
\end{equation}
It is well known that these bilinear forms are continuous and coercive with respect to the $H^1$-norm, the latter being dependent on Korn's first inequality \cite{ciarlet}:
\begin{equation} \label{eq:continuity_bilinear}
\begin{aligned}
&\exists \  m_A>0, m_B>0 \ \text{such that} \ \forall \ \bm{u},\bm{v} \in \bm{V}, \\
    &\lvert A[\bm{u},\bm{v}] \rvert \leq m_A\lVert \bm{u} \rVert_{\bhonepm}\lVert \bm{v} \rVert_{\bhonepm}, \ \lvert B[\bm{u},\bm{v}]\rvert \leq m_B\lVert \bm{u} \rVert_{\bhonepm}\lVert \bm{v} \rVert_{\bhonepm},
\end{aligned}
\end{equation}
\begin{equation} \label{eq:v-ellipticity}
    \exists \  \alpha_A>0,\alpha_B>0 \ \text{such that} \ \forall \ \bm{v}\in \bm{V}, \   A[\bm{v},\bm{v}]\geq\alpha_A\lVert \bm{v}\rVert^2_{\bhonepm}, \ B[\bm{v},\bm{v}]\geq \alpha_B\lVert\bm{v}\rVert^2_{\bhonepm}. 
\end{equation}
We also define two nonlinear functionals, $J_{\nu}[\cdot,\cdot]: \bm{V} \times \bm{V} \to \mathbb{R}$ and $J_{\tau}[\cdot,\cdot]: \bm{V} \times \bm{V} \to \mathbb{R}$, representing the contributions from the normal compliance law and the Coulomb friction law, respectively:
\begin{equation} \label{eq:j_norm_functional}
    J_{\nu}[\bm{u},\bm{v}]=\int_{\mc}c_{\nu}\left([\bu]_{\nu}-g\right)_+[\bv]_{\nu} ds,
\end{equation}
\begin{equation} \label{eq:j_tang_functional}
    J_{\tau}[\bm{u},\bm{v}]=\int_{\mc}\mu_{\text{fr}} c_{\nu}([\bu]_{\nu}-g)_+\lvert[\bm{v}]_{\tau}\rvert ds,
\end{equation}
and we collect these into a single nonlinear functional $J[\cdot,\cdot]: \bm{V} \times \bm{V} \to \mathbb{R}$:
\begin{equation} \label{eq:j_functional}
    J[\bm{u},\bm{v}]=J_{\nu}[\bu,\bv]+J_{\tau}[\bu,\bv].
\end{equation}
The variational problem, henceforth referred to as problem (Q), then reads as follows: 

\begin{problem}
Given $\bm{f} \in L^2(0,T;\bm{L}^2(\Omega))$, $\psi \in L^2(0,T;L^2(\Omega))$, $\sfrac \in L^2(0,T;L^2(\mathcal{C}))$ and $g \in L^{\infty}(\mc)$, find $\bm{u} \in H^1(0,T;\bm{V})$, $p \in L^2(0,T;Q) \cap H^1(0,T;Q^{\prime})$ and $\pc \in   L^2(0,T;H^1(\mathcal{C})) \cap H^1(0,T;H^{-1}(\mc))$, such that, for a.e. $t \in (0,T)$:
\begin{equation} \label{eq:weak_mechanics}
\begin{aligned}
    A[\bu,\bv-\dotu]+B[\dotu,\bv-\dotu] 
    &-\alpha (p,\nabla \cdot (\bm{v}-\dotu))\\
    &+(\pc, [\bm{v}-\dotu]_{\nu})_{\mathcal{C}}
    +J[\bu,\bv]-J[\bu,\dotu] \geq (\bm{f},\bm{v}-\dotu) \ \ \forall\bm{v} \in \bm{V},
\end{aligned}
\end{equation}
\begin{equation} \label{eq:weak_biot}
\begin{aligned}
    \langle c_0\dotp,q \rangle +(\alpha\nabla \cdot \dotu,q)+\frac{1}{\mu}(\bk\nabla p, \nabla q)=(\psi,q)+\frac{1}{6\mu}\sum_{*\in \pm}(g(\pc-p^*),q^*)_{\mc} \ \ \forall q \in Q,
\end{aligned}
\end{equation}
\begin{equation} \label{eq:weak_biot_fracture}
    \begin{aligned}
        \langle c_{f}g\dot{\pi},\theta \rangle_{\mc}-([\dotu]_{\nu}, \qc )_{\mathcal{C}}&+\frac{1}{12\mu}(g^3 \tnabla \pc,\tnabla \qc )_{\mathcal{C}}\\
        &=( g\sfrac,\qc)_{\mc}-\frac{1}{6\mu}\sum_{*\in\pm}\left(g(\pc-p^*),\qc\right)_{\mc} \ \ \forall \qc \in H^1(\mc).
    \end{aligned}
\end{equation}
subject to the initial conditions:
\begin{align}
    \bu\lvert_{t=0} &= \bu_0, \label{eq:init_cond1}\\
    (c_0p+\alpha\nabla \cdot \bu)\lvert_{t=0} &= c_0p_0+\alpha \nabla \cdot \bu_0, \label{eq:init_cond2}\\
    (c_{f}g\pc-[\bu]_{\nu})\rvert_{t=0}&=c_{f}g\pcnull-[\bu_0]_{\nu} \label{eq:init_cond3}.
\end{align}
\end{problem}
Note that the modeling of a viscoelastic solid allows the displacement to be taken in a stronger temporal space, compared to a purely elastic solid.

The variational mass balance equation \eqref{eq:weak_biot} is derived by assuming a sufficiently regular solution to the equation in strong form \eqref{eq:biot_mass_balance}, multiplying by $q\in Q$, integrating over and applying Green's formula separately in $\Omega^+$ and $\Omega^-$. The resulting equation obtained after adding the equations in $\Omega^+$ and $\Omega^-$ is the following:
\begin{equation}
    (c_0\dot{p}+\alpha\nabla \cdot \dot{\bm{u}},q)+\frac{1}{\mu}(\bk \nabla p,\nabla q)+\sum_{*\in \pm}(\bv_D^*\cdot \bm{\nu}^*,q^*)_{\mc}=(\psi,q),
\end{equation}
where the boundary integrals reduce to integrals over $\mc$ (the third term on the left-hand side), due to the homogeneous boundary conditions \eqref{eq:bcs} and the continuity conditions \eqref{eq:interface_cond_flow} of pressure and normal flux across $\Gamma\setminus\mc$.  Recalling the interface condition \eqref{eq:normal_darcy} for the normal flux on $\mc$, equation \eqref{eq:weak_biot} follows. The variational fracture mass balance equation \eqref{eq:weak_biot_fracture} is derived in a similar manner, by multiplying \eqref{eq:mass_balance_fracture} by $\qc\in H^1(\mc)$, integrating over $\mc$ and using Green's formula on the diffusive term. The boundary term vanishes due to the homogeneous Neumann condition on the fracture tips \eqref{eq:fracture_tips}.

Next, we explain the derivation of the variational inequality \eqref{eq:weak_mechanics}. By multiplying the momentum balance equation \eqref{eq:mom_bal} with $\bm{v}-\dotu, \bm{v} \in \bm{V}$, integrating and applying Green's formula separately in $\Omega^+$ and $\Omega^-$, we get, after adding the equations on each subdomain:
\begin{equation} \label{eq:weak_elasticity_derivation}
\begin{aligned}
    A[\bu,\bv-\dotu]+B[\dotu,\bv-\dotu]-\alpha( p,\nabla \cdot (\bm{v}-\dotu))=(\bm{f},\bm{v}-\dotu)+\sum_{*\in \pm}\big((\bm{\sigma}^{\mathrm{por}})^*\cdot\bm{\nu}^*,(\bm{v}-\dotu)^*\big)_{\mathcal{C}}.
\end{aligned}
\end{equation}
Once again, the boundary integrals reduce to integrals over $\mc$, due to the homogeneous boundary conditions \eqref{eq:bcs} and the continuity conditions \eqref{eq:interface_cond_mech} of displacement and normal stress across $\Gamma\setminus\mc$. Additionally, since equations \eqref{eq:contact_traction} and \eqref{eq:interface_force_balance} imply that $[\bm{\sigma}^{\mathrm{por}}]_{\nu}=\bm0$ over $\mc$, the boundary integrals can be rewritten as
\begin{equation}
    \sum_{*\in \pm}\big((\bm{\sigma}^{\mathrm{por}})^*\cdot\bm{\nu}^*,(\bm{v}-\dotu)^*\big)_{\mathcal{C}}=((\bm{\sigma}^{\text{por}})^+ \cdot \bm{\nu}^+,[\bv-\dotu])_{\mc}.
\end{equation}
Using \eqref{eq:contact_traction} and splitting up into normal and tangential components, this is further rewritten as
\begin{equation} \label{eq:boundary_integral}
\begin{aligned}    \left((\bm{\sigma}^{\mathrm{por}})^+\cdot\bm{\nu}^+,[\bm{v}-\dotu]\right)_{\mc}&=(\bm{t}^+,[\bv-\dotu])_{\mc}-(\pc\bm{\nu}^+,[\bv-\dotu])_{\mc}\\
&=(t_{\nu},[\bv-\dotu]_{\nu})_{\mc}+(\bm{t}_{\tau},[\bv-\dotu]_{\tau})_{\mc}-(\pc,[\bv-\dotu]_{\nu})_{\mc},
\end{aligned}
\end{equation}
Next, the Coulomb friction law \eqref{eq:friction_law} implies the inequality:
\begin{equation} \label{eq:friction_inequality}
    \bm{t}_{\tau} \cdot [\bv-\dotu]_{\tau}\geq \mu_{\text{fr}} c_{\nu}([\bu]_{\nu}-g)_+(\lvert [\dotu]_{\tau}\rvert-\lvert[\bv]_{\tau}\rvert),
\end{equation}
see e.g. chapter 10 of \cite{kikuchi_oden} for a proof. From \eqref{eq:boundary_integral}, \eqref{eq:friction_inequality} and the normal compliance law \eqref{eq:normal_compliance} we obtain
\begin{equation} \label{eq:friction_inequality_2}
    \left((\bm{\sigma}^{\mathrm{por}})^+\cdot\bm{\nu}^+,[\bm{v}-\dotu]\right)_{\mc}\geq J[\bu,\dotu]-J[\bu,\bv]-(\pc,[\bv-\dotu]_{\nu})_{\mc},
\end{equation}
using the previously defined $J$-functional \eqref{eq:j_functional}.
Combining \eqref{eq:friction_inequality_2} and \eqref{eq:weak_elasticity_derivation} yields \eqref{eq:weak_mechanics}.

Regarding the question of existence and uniqueness of solutions to problem (Q), we note that existence and uniqueness for a very similar problem is proven in the paper of de Hoop and Kumar \cite{dehoop2023}. They consider a mixed-dimensional model of poroviscoelasticity including fracture contact mechanics, using the same contact conditions as in the present paper. The key step in their proof of existence is the fixed point iteration between the flow and mechanics equations. The two differences between their model and ours is that they consider the fully dynamic Biot equations, with an acceleration term included in the momentum balance equation, and they assume continuity of the pressure across the fracture. Both of these differences, however, do not interfere with the arguments presented in the cited reference. The pressure conditions across the fracture-matrix interface in our case leads to a positive contribution on the left hand side, as shown in the estimates of Theorems \ref{thm:stability_estimate} and \ref{thm:error_estimate}. Further, the lack of acceleration term in the momentum balance equation does not cause any difficulties in the fixed point iteration scheme, as argued in our earlier work \cite{nevland2026numerical}. We thus expect that the proof can be adapted to our setting. However, the focus of this work is on the numerical analysis and we simply assume the the well-posedness of the continuous formulation. 


\section{Discretization} \label{sec:discretization}

In this section, we present a fully discrete scheme for numerically approximating solutions to problem (Q) of Section \ref{sec:variational_form}, using conforming finite elements in space and the implicit Euler method in time. We choose $N$ equidistant time points in the interval $[0,T]$. The time step size is accordingly $\Delta t=\frac{T}{N}$ and the discrete time points are $t_n=n\Delta t, \ 0 \leq n \leq N$. For any time-dependent function $v(\cdot,t)$, assumed to be continuous in time, we write $v_n:=v(\cdot,t_n)$. The backward finite difference approximation of the time derivative is denoted as follows:
\begin{equation} \label{eq:time_difference}
    \delta_t v_n :=\frac{v_n-v_{n-1}}{\Delta t}.
\end{equation}
Let $\mathcal{T}^h$ denote a regular family of conforming triangulations of $\overline{\Omega}$, indexed by the maximal diameter $h$ of the elements of $\mathcal{T}^h$. For simplicity, we assume that $\mathcal{T}^h$ meshes $\Omega^+$ and $\Omega^-$ separately, implying that $\mathcal{C}$ does not cross any elements of $\mathcal{T}^h$. Moreover, we assume that the fracture $\mc$ is meshed by the trace of $\mathcal{T}^h$ on $\mc$, call this mesh $\mc^h$. Let $k_1\geq1$, $k_2\geq1$ and $k_3 \geq 1$ denote polynomial orders. We consider standard Lagrangian finite elements, where the displacement space $\bm{V}$, the matrix pressure space $Q$ and the fracture pressure space $H^1(\mc)$ are approximated by $\bm{V}^h$, $Q^h$ and $\Pi^h$, respectively, defined as follows:
\begin{align} \label{eq:fem_spaces}
    \bm{V}^h&=\{ \bv^h \in \bm{V} : \forall \ T\in\mathcal{T}^h, \ \bv^h\rvert_{T}\in \mathbb{P}^d_{k_1}(T) \},\\
    Q^h &= \{ q^h\in Q : \forall \ T\in\mathcal{T}^h, \ q^h\rvert_{T}\in \mathbb{P}_{k_2}(T) \},\\
    \Pi^h &= \{\qch \in H^1(\mc) : \forall e \in \mc^h, \ \qch\lvert_e \in \mathbb{P}_{k_3}(e) \}.
\end{align}
where $\mathbb{P}_p(T)$ denotes the space of polynomials of order $p$ on $T$ ($\mathbb{P}_p^d(T)$ denotes the Cartesian product of $d$ such spaces). We also introduce suitable approximation operators $I_{\bm{V}^h}:\bm{V}\to \bm{V}^h$, $I_{Q^h}: Q\to Q^h$ and $I_{\Pi^h}:H^1(\mc)\to \Pi^h$, and use the notation $\bu^I := I_{\bm{V}^h}(\bu)$, $p^I:=I_{Q^h}(p)$, $\pc^I := I_{\Pi^h}(\pc)$. These operators are assumed to satisfy the following error estimates on $\Omega^{\pm}$:
\begin{equation}
\begin{aligned}
\label{eq:interpolation_estimate}
    &\lVert\bm{u}-\bu^I \rVert_{\bhonepm}\leq C_1h^{q-1}\lvert\bu \rvert_{\bm{H}^{q}(\Omega^{\pm})}, \ \lVert p-p^I \rVert_{\honepm} \leq C_2h^{r-1}\lvert p \rvert_{H^{r}(\Omega^{\pm})},\\
    &\lVert \pc-\pc^I \rVert_{H^1(\mc)} \leq C_3h^{s-1}\lvert \pc \rvert_{H^{s}(\mc)},
\end{aligned}
\end{equation}
assuming that $\bu \in \bm{H}^{q}(\Omega)$, $p \in H^{r}(\Omega)$ and $\pc \in H^{s}(\mc)$, where
\begin{equation} \label{eq:sobolev_exponents}
    q \in (1,k_1+1], \ r \in (1, k_2+1], \ s \in(1, k_3+1].
\end{equation}
Note that $q, r$ and $s$ are not necessarily integers, and so the spaces for $\bu,p$ and $\pc$ may be fractional Sobolev spaces. An example of an operator satisfying the above estimates is the Scott-Zhang interpolation operator \cite{scott_zhang}.

\subsection{Fully discrete approximation} \label{sec:fully_discrete_form}

The following fully discrete conformal variational problem is proposed, called problem (Q-D).

\begin{discreteproblem}
Assume $\bm{f}\in H^1(0,T;\bm{L}^2(\Omega))$, $\psi\in H^1(0,T;L^2(\Omega))$, $\sfrac\in H^1(0,T;L^2(\mc))$ and $g\in L^{\infty}(\mc)$.\\

At time $t=0$, given initial data $\bu_0 \in \bm{V}$, $p_0 \in Q$ and $\pcnull\in H^1(\mc)$, set $\bu_0^h =\bu_0^I$, $p_0^h=p_0^I$ and $\pcnull^{h}=\pcnull^I$. For any $n \in \{1,2,..,N\}$, given $\bm{u}_{n-1}^h, p_{n-1}^h$ and $\pcnp$, find $\bm{u}_n^h \in \bm{V}^h$, $p_n^h \in Q^h$ and $\pch \in \Pi^h$ such that:
\begin{equation} \label{eq:discrete_mechanics}
\begin{aligned}
    A[\bm{u}_n^h,\bm{v}^h&-\dunh]+B[\dunh,\bm{v}^h-\dunh]-\alpha (p_n^h,\nabla \cdot (\bm{v}^h-\dunh))\\
    &+(\pch,[\bv^h-\dunh]_{\nu})_{\mc}+J[\bm{u}_n^h,\bm{v}^h]-J[\bm{u}_n^h,\dunh] \geq ( \bm{f}_n,\bm{v}^h-\dunh) \ \ \forall \bm{v}^h \in \bm{V}^h,
\end{aligned}
\end{equation}
\begin{equation} \label{eq:discrete_biot}
\begin{aligned}
    \left( c_0\dpnh+\alpha\nabla \cdot \dunh, q^h\right)&+\frac{1}{\mu}(\bk\nabla p_n^h,\nabla q^h)\\
    &=(\psi_n,q^h)+\frac{1}{6\mu}\sum_{*\in\pm}\big(g(\pch-\pnhtr),\qhtr\big)_{\mc} \ \ \forall q^h \in Q^h,
\end{aligned}
\end{equation}
\begin{equation} \label{eq:discrete_biot_fracture}
\begin{aligned}
    (c_fg\delta_t \pch-[\dunh]_{\nu},\qch)_{\mc}&+\frac{1}{12\mu}(g^3\tnabla \pch,\tnabla \qch)_{\mc}\\
    &=( g\sfrac_n,\qch)_{\mc}-\frac{1}{6\mu}\sum_{*\in\pm}\big(g(\pch-\pnhtr),\qch\big)_{\mc} \ \ \forall \qch \in \Pi^h.
\end{aligned}
\end{equation}
\end{discreteproblem}
The assumed regularity in time of the data ensures that $\bm{f}_n$, $\psi_n$ and $\sfrac_n$ are well-defined.

\subsection{Existence and uniqueness of the discrete solution}

To prove the existence and uniqueness of the solution to problem (Q-D), we follow the same strategy as our previous paper \cite{nevland2026numerical}. In that paper, the existence and uniqueness of the Biot equations coupled to frictional contact mechanics (in a non-fractured medium) was proven by first considering a simpler model where the contact mechanics terms are lagged. Existence and uniqueness of this simpler problem was proven with help of the fixed-stress splitting scheme, which is a well-known algorithm for solving problems of coupled flow and mechanics, by decoupling the flow and mechanics parts and iterating between the two problems \cite{kim2011stability,mikelic2013convergence}. Finally, the existence and uniqueness of the full problem including Coulomb friction was proven by a fixed-point argument. The extension of these techniques to the mixed-dimensional setting is relatively straightforward, with the most significant changes occurring in the contraction proof of the fixed-stress splitting scheme, which becomes more involved when fracture flow is included \cite{girault_splitting}.

Our simplified model consists of replacing the variational inequality \eqref{eq:discrete_mechanics} by
\begin{equation} \label{eq:simple_var_ineq}
\begin{aligned}
    A[\bm{u}_n^h,\bm{v}^h&-\dunh]+B[\dunh,\bm{v}^h-\dunh]-\alpha (p_n^h,\nabla \cdot (\bm{v}^h-\dunh))\\
    &+(\pch,[\bv^h-\dunh]_{\nu})_{\mc}+J[\bm{\eta},\bm{v}^h]-J[\bm{\eta},\dunh] \geq ( \bm{f}_n,\bm{v}^h-\dunh) \ \ \forall \bm{v}^h \in \bm{V}^h,
\end{aligned}
\end{equation}
for a given $\bm{\eta}\in\bm{V}^h$. That is, the first arguments of the $J$-functionals are kept fixed in this simplified variational inequality.  Recalling the definitions of the $J$-functional \eqref{eq:j_functional} and the normal compliance law \eqref{eq:normal_compliance}, we observe that this simplification corresponds to prescribing the normal contact traction on the fracture surface, leading to the so-called Tresca friction model, which is known to be feasible for analysis \cite{chouly2023finite,kikuchi_oden}. The function $\bm{\eta}$ may be chosen differently for different time steps, but for simplicity, we omit the subscript $n$. The full simplified problem can be summarized as follows:
\begin{simpleproblem}
    Keep the assumptions and initial conditions of problem (Q-D). For any $n\in \{1,2,...,N\}$, given $\bu_{n-1}^h, p_{n-1}^h, \pcnp$ and $\bm{\eta}$, find $\bu_n^h\in\bm{V}^h$, $p_n^h\in Q^h$ and $\pch\in \Pi^h$ solving \eqref{eq:discrete_biot}, \eqref{eq:discrete_biot_fracture} and \eqref{eq:simple_var_ineq}.
\end{simpleproblem}
We use the subscript $k$ to denote the iteration of the splitting scheme. Moreover, we reuse the $\delta_t$-notation in the splitting scheme, with the meaning 
\begin{equation}
   \delta_t v_{n,k}^h:=\frac{v_{n,k}^h-v_{n-1}^h}{\Delta t},
\end{equation}
for any function $v_{n,k}$ depending on the time step $t_n$, the mesh $\mathcal{T}^h$ and the splitting scheme iteration $k$.
The fixed-stress splitting scheme for problem (Q-D-S) reads:

\textbf{Step 1: Flow}. Given $\ppfs\in Q^h$, $\pcpfs \in \Pi^h$ and $\upfs \in \bm{V}^h$, find $\pfs \in Q^h$ and $\pcfs \in \Pi^h$ such that
\begin{equation} \label{eq:fixed_stress_biot}
\begin{aligned}
    \forall q^h \in Q^h \ , \ &\left( (c_0+\frac{\alpha^2}{\lambda})\delta_t\pfs+\alpha\nabla \cdot \delta_t\upfs, q^h\right)+\frac{1}{\mu}\left(\bk\nabla \pfs,\nabla q^h\right)\\
    &=\left(\frac{\alpha^2}{\lambda}\delta_t \ppfs,q^h\right)+\left(\psi_n,q^h\right) + \frac{1}{6\mu}\sum_{*\in\pm}\left(g(\pcfs-\pfstr),\qhtr\right)_{\mc},
\end{aligned}
\end{equation}
\begin{equation} \label{eq:fixed_stress_biot_fracture}
    \begin{aligned}
        \forall \qch \in \Pi^h, \ &\big((c_{f}g+\gamma_c)\delta_t\pcfs-[\delta_t\upfs]_{\nu},\qch\big)_{\mc}+\frac{1}{12\mu}\left(g^3 \tnabla \pcfs,\tnabla \qch\right)_{\mc}\\
        &=\left(\gamma_c \delta_t \pcpfs,\qch\right)_{\mc}+\left( g \sfrac_n,\qch\right)_{\mc}- \frac{1}{6\mu}\sum_{*\in\pm}\left(g(\pcfs-\pfstr),\qch\right)_{\mc}
    \end{aligned}
\end{equation}
\textbf{Step 2: Mechanics}. Given $\pfs,\pcfs$ from step 1, find $\ufs \in \bm{V}^h$ such that
\begin{equation} \label{eq:fixed_stress_mechanics}
\begin{aligned}
    \forall \bm{v}^h \in \bm{V}^h \ , \ &A[\ufs,\bm{v}^h-\delta_t \ufs]+B[\delta_t \ufs,\bm{v}^h-\delta_t \ufs]\\
    &-\alpha (\pfs,\nabla \cdot (\bm{v}^h-\delta_t \ufs))+(\pcfs,[\bv^h-\delta_t \ufs]_{\nu})_{\mc}\\
    &+J[\bm{\eta},\bm{v}^h]-J[\bm{\eta},\delta_t \ufs]\geq ( \bm{f}_n,\bm{v}^h-\delta_t \ufs),
\end{aligned}
\end{equation}
and the scheme is initialized for time step $n \geq 1$ by setting $p_{n,0}^{h}=p_{n-1}^h$, $\pi_{n,0}^{h}=\pi_{n-1}^h$ and $\bm{u}_{n,0}^{h}=\bm{u}_{n-1}^h$. Note that in equations \eqref{eq:fixed_stress_biot} and \eqref{eq:fixed_stress_biot_fracture}, we have added the stabilization terms $\alpha^2\lambda^{-1}\delta_t \pfs$ and $\gamma_c\delta_t \pcfs$ to the left-hand side, with similar terms on the right-hand side for consistency. The first term is a standard addition for the fixed-stress splitting scheme for the Biot equations \cite{mikelic2013convergence}, and the second term was introduced in \cite{girault_splitting} as part of an extension of the splitting scheme to the Biot equations in fractured porous media. A formula for the stabilization parameter $\gamma_c$, depending on material parameters and the minimal fracture gap $g_{\text{min}}$, is established in the upcoming contraction proof of the fixed-stress scheme.

It is straightforward to check that each step of the above fixed-stress splitting scheme has a unique solution. The first step, \eqref{eq:fixed_stress_biot} and \eqref{eq:fixed_stress_biot_fracture}, corresponds to a square system of linear equations in finite dimension, from which the existence and uniqueness follows by checking the invertibility of the corresponding matrix (we skip the proof as this is a straightforward computation; see e.g. our previous paper \cite{nevland2026numerical}). The second step \eqref{eq:fixed_stress_mechanics} corresponds to a standard variational inequality of the second kind, since the pressures are fixed during this step, and the existence and uniqueness follows from standard results on variational inequalities; see for instance \cite{han2002quasistatic}.

The following Lemma establishes the existence and uniqueness of a solution to problem (Q-D-S), which coincides with the limit functions of the fixed-stress splitting scheme. It is proven in Appendix \ref{app:contraction}, with the help of a contraction estimate on the fixed-stress iterations.
\begin{lemma} \label{lemma:fixed_stress_limits}
    Given a sufficiently small time step $\Delta t$, there exists $p_n^h \in Q^h$, $\pch \in \Pi^h$ and $\bm{u}_n^h \in \bm{V}^h$ such that the iterations $\ppfs, \pcpfs, \upfs$ of the fixed-stress splitting scheme satisfy, as $k\to\infty$: 
    \begin{equation*}
        \ppfs \to p_n^h \ \text{in} \ \honepm, \ \ 
        \pcpfs \to \pch \ \text{in} \ H^1(\mc), \ \ 
        \upfs\to \bm{u}_n^h \ \text{in} \ \bhonepm,
    \end{equation*}
    with strong convergence in the above norms. Moreover, the limit functions $\bu_n^h, p_n^h$ and $\pch$ solve problem (Q-D-S) at time step $n$.
\end{lemma}
We will also need the following Lemma regarding an estimate on the $J$-functionals. The proof is skipped, as it is a straightforward extension of a similar estimate in our previous paper \cite{nevland2026numerical} to the mixed-dimensional setting.
\begin{lemma} \label{lemma:friction_estimate}
    The $J$-functional \eqref{eq:j_functional} satisfies the following estimate:
    \begin{equation}
        J[\bu_1,\bv_2]-J[\bu_1,\bv_1]+J[\bu_2,\bv_1]-J[\bu_2,\bv_2] \leq C_J \lVert \bu_1-\bu_2 \rVert_{\bhonepm}\lVert \bv_1-\bv_2 \rVert_{\bhonepm},
    \end{equation}
    for all $\bu_1,\bu_2,\bv_1,\bv_2 \in \bm{V}$, for some $C_J>0$ depending on $c_{\nu}, \mu_{\textnormal{fr}}$ and the constants $C^+, C^-$ of the trace inequality \eqref{eq:trace_inequality}.
\end{lemma}
With the help of the two preceding Lemmas, we can prove the existence and uniqueness of a solution to the original problem (Q-D).
\begin{theorem} \label{thm:existence_discrete}
    For a sufficiently small time step $\Delta t$, problem (Q-D) is guaranteed to have a unique solution. In other words, there exists for any $n\in \{1,2,...,N\}$ a unique tuple $(\bu_n^h,p_n^h, \pch)\in \bm{V}^h \times Q^h \times \Pi^h$ that solves \eqref{eq:discrete_mechanics}, \eqref{eq:discrete_biot} and \eqref{eq:discrete_biot_fracture}.
\end{theorem}
\begin{proof}
The existence of a solution is proven by performing fixed-point iteration on the given function $\bm{\eta}$ in problem (Q-D-S); this fixed-point argument is taken from an article \cite{Han1999} concerning a purely mechanical contact problem. Introduce the operator $\Lambda:\bm{V}^h \to \bm{V}^h$, defined by
\begin{equation} \label{eq:fixed_point_operator}
    \Lambda(\bm{\eta})=\bm{u}^{\bm{\eta}}_n \ \ \ \forall\bm{\eta}\in \bm{V}^h,
\end{equation}
where $\bm{u}^{\bm{\eta}}_n$ denotes the displacement solution to problem (Q-D-S) for a given function $\bm{\eta}$ (i.e. the limit of the sequence $\upfs$ of the fixed-stress splitting scheme as $k\to\infty$, as obtained in Lemma \ref{lemma:fixed_stress_limits}). We similarly denote the corresponding pressure solutions by $p_n^{\bm{\eta}}$ and $\pi_{n}^{\bm{\eta}}$. For ease of notation, we drop the $h$ superscript in this proof. We will show that this operator has a unique fixed point.

Consider two solutions $(\bm{u}^{\bm{\eta_1}}_n, p^{\bm{\eta_1}}_n, \pi_{n}^{\bm{\eta_1}})$ and $(\bm{u}^{\bm{\eta_2}}_n, p^{\bm{\eta_2}}_n, \pi_{n}^{\bm{\eta_2}})$ to problem (Q-D-S), corresponding to different functions $\bm{\eta}_1$ and $\bm{\eta}_2$, respectively. We use a tilde to denote the differences between these solutions; $\widetilde{p}_n=p_n^{\bm{\eta}_1}-p_n^{\bm{\eta}_2}$, $\widetilde{\pi}_{n}=\pi_{n}^{\bm{\eta}_1}-\pi_{n}^{\bm{\eta}_2}$, $\widetilde{\bm{u}}_n=\bm{u}^{\bm{\eta}_1}_n-\bm{u}^{\bm{\eta}_2}_n$ and $\widetilde{\bm{\eta}}=\bm{\eta}_1-\bm{\eta}_2$. The simplified variational inequality \eqref{eq:simple_var_ineq} corresponding to the $\bm{u}^{\bm{\eta}_1}_n$-solution is tested with $\bm{v}^h=\delta_t\bm{u}^{\bm{\eta}_2}_n$, and the inequality corresponding to the $\bm{u}^{\bm{\eta}_2}_n$-solution is tested with $\bm{v}^h= \delta_t\bm{u}^{\bm{\eta}_1}_n$, and the resulting equations are added. The discrete mass balance equations \eqref{eq:discrete_biot} and \eqref{eq:discrete_biot_fracture} are tested with $q^h=\widetilde{p}_n$ and $\qch=\widetilde{\pi}_{n}$, and the equations for each solution are subtracted. After adding the flow and mechanics equations, the resulting estimate becomes:
\begin{equation} \label{eq:fpi_arg_1}
\begin{aligned}
    &A[\widetilde{\bm{u}}_n, \widetilde{\bm{u}}_n]+\frac{1}{\Delta t}B[\widetilde{\bm{u}}_n, \widetilde{\bm{u}}_n]+c_0\lVert \widetilde{p}_n \rVert_{\ltwo}^2+c_{f}\lVert g^{1/2}\widetilde{\pi}_{n}\rVert_{L^2(\mc)}^2+\frac{\Delta t}{\mu}\lVert \bm{K}^{1/2}\nabla\widetilde{p}_n\rVert_{\bltwopm}^2\\
    &+\frac{\Delta t}{12\mu}\lVert g^{3/2}\tnabla \widetilde{\pi}_{n}\rVert_{\bltwoc}^2+\frac{\Delta t}{6\mu}\sum_{*\in\pm} \lVert g^{1/2}(\widetilde{\pi}_{n}-\widetilde{p}_n^*)\rVert_{L^2(\mc)}^2\\
    &\leq J[\bm{\eta}_1,\bm{u}^{\bm{\eta}_2}_n-\bm{u}_{n-1}]
    -J[\bm{\eta}_1,\bm{u}^{\bm{\eta}_1}_n-\bm{u}_{n-1}]
    +J[\bm{\eta}_2,\bm{u}^{\bm{\eta}_1}_n-\bm{u}_{n-1}]
    -J[\bm{\eta}_2,\bm{u}^{\bm{\eta}_2}_n-\bm{u}_{n-1}]\\
    &\leq C_J\lVert \widetilde{\bm{\eta}}\rVert_{\bhonepm}\lVert \widetilde{\bm{u}}_n\rVert_{\bhonepm},
\end{aligned}
\end{equation}
where Lemma \ref{lemma:friction_estimate} has been used for the final inequality. Together with the coercivity conditions \eqref{eq:v-ellipticity}, this implies the estimate
\begin{equation} \label{eq:fpi_arg_3}
    (\alpha_A+\frac{\alpha_B}{\Delta t})\lVert \widetilde{\bm{u}}_n\rVert_{\bhonepm} \leq C_J\lVert \widetilde{\bm{\eta}}\rVert_{\bhonepm},
\end{equation}
From this we obtain that the operator $\Lambda$ \eqref{eq:fixed_point_operator} is a contraction, and consequently has a unique fixed point, if
\begin{equation} \label{eq:time_step_constraint_1}
   \alpha_A+\frac{\alpha_B}{\Delta t}>C_J,
\end{equation}
which can be achieved by choosing a sufficiently small time step. Denote this fixed point by $\overline{\bm{\eta}}$, and let $(\overline{\bm{u}}_n, \overline{p}_n, \overline{\pi}_{n})$ denote the corresponding solution to problem (Q-D-S) with $\overline{\bm{\eta}}$ as the given function. Then $(\overline{\bm{u}}_n, \overline{p}_n, \overline{\pi}_{n})$ is also a solution to problem (Q-D). Indeed, the only difference between problem (Q-D) and problem (Q-D-S) is the variational inequality, but due to the fixed-point property $\overline{\bm{\eta}}=\overline{\bm{u}}_n$, the solution to the simplified variational inequality \eqref{eq:simple_var_ineq} coincides with that of the full variational inequality \eqref{eq:discrete_mechanics}. We thus conclude that there exists a solution to problem (Q-D), given a sufficiently small time step.

Finally, we prove the uniqueness of the solution to problem (Q-D). This proof is straightforward at this point, as it follows the exact same steps as the preceding proof of the contraction of the operator $\Lambda$. Assuming there are two solutions, $(\bm{u}_{n1}^h,p_{n1}^h, \pi_{n1}^h)$ and $(\bm{u}_{n2}^h,p_{n2}^h, \pi_{n2}^h)$, to the discrete equations \eqref{eq:discrete_mechanics}, \eqref{eq:discrete_biot} and \eqref{eq:discrete_biot_fracture}, denoting the differences between the solutions by tildes, and performing the same computations as before, \eqref{eq:fpi_arg_1}-\eqref{eq:fpi_arg_3}, we obtain the estimate
\begin{equation} \label{eq:uniqueness_estimate}
\begin{aligned}
    &\left(\alpha_A+\frac{\alpha_B}{\Delta t}\right)\lVert \tilu\rVert_{\bhonepm}^2+c_0\lVert \tilp \rVert_{\ltwo}^2+c_{f}\lVert g^{1/2}\widetilde{\pi}_{n}^h\rVert_{L^2(\mc)}^2+\frac{\Delta t}{\mu}\lVert \bk^{1/2}\nabla \tilp \rVert_{\bltwopm}^2\\
    &+\frac{\Delta t}{12\mu}\lVert g^{3/2}\tnabla \widetilde{\pi}_{n}^h\rVert_{\bltwoc}+\frac{\Delta t}{6\mu}\sum_{*\in\pm} \lVert g^{1/2}(\widetilde{\pi}_{n}^h-\widetilde{p}_n^{h*})\rVert_{L^2(\mc)}^2 \leq C_J\lVert \tilu \rVert_{\bhonepm}^2,
\end{aligned}
\end{equation}
from which the uniqueness follows, given the fulfillment of the same condition as before, \eqref{eq:time_step_constraint_1}. Theorem \ref{thm:existence_discrete} follows.
\end{proof}

\subsection{Stability estimate} \label{sec:stability_estimate}

In this section, we prove the stability of the solution $(\bu_n^h, p_n^h, \pch)$ to the discrete problem (Q-D).

\begin{theorem} \label{thm:stability_estimate}
   If the time step size satisfies $\Delta t<2\alpha_A\alpha_B$, then the solution to problem (Q-D) satisfies the following stability bound:
\begin{equation} \label{eq:stability_estimate}
\begin{aligned}
    \lVert \bu_N^h\rVert_{\bhonepm}^2&+\lVert p_N^h \rVert_{\ltwo}^2+\lVert \pcnh \rVert_{L^2(\mc)}^2+\sum_{n=1}^N \Delta t \lVert \delta_t \bu_n^h \rVert_{\bhonepm}^2\\
   & + \sum_{n=1}^N \Delta t \lVert \nabla p_n^h \rVert_{\bltwopm}^2+\sum_{n=1}^N \Delta t \lVert \tnabla \pch \rVert_{\bltwoc}^2\leq C,
\end{aligned}
\end{equation}
where $C$ depends on the data and the final time $T$, but is independent of $n, h$ and $\Delta t$.
\end{theorem}
\begin{proof}
    We shall need the following discrete version of the chain rule, valid for some arbitrary $v_n$ defined at the discrete time points $t_n$:
\begin{equation} \label{eq:discrete_chain_rule}
    (g^{1/2}\delta_t v_n,g^{1/2}v_n)=\frac{1}{2\Delta t}\left( \lVert g^{1/2}v_n \rVert_{\ltwopm}^2-\lVert g^{1/2}v_{n-1}\rVert_{\ltwopm}^2+\lVert g^{1/2}(v_n-v_{n-1})\rVert_{\ltwopm}^2\right).
\end{equation}
The flow equation in the matrix \eqref{eq:discrete_biot} is tested with $q^h=p_n^h$. Using the above chain rule (with $g=1$), the obtain the following stability equality:
\begin{equation} \label{eq:stability_flow_1}
    \begin{aligned}
        &\frac{c_0}{2\Delta t}\left(\lVert p_n^h \rVert_{\ltwo}^2-\lVert p_{n-1}^h \rVert_{\ltwo}^2+\lVert  p_n^h-p_{n-1}^h \rVert_{\ltwo}^2\right)+\alpha\left(p_n^h,\nabla \cdot \frunh\right)\\
        &+\frac{1}{\mu}\lVert \bk^{1/2}\nabla p_n^h \rVert_{\bltwopm}^2=(\psi_n,p_n^h)+\frac{1}{6\mu}\sum_{*\in\pm}\Big(g^{1/2}(\pch-\pnhtr),g^{1/2}\pnhtr\Big)_{\mc}.
    \end{aligned}
\end{equation}
The flow equation in the fracture \eqref{eq:discrete_biot_fracture} is similarly tested with $\qch=\pch$, and the chain rule \eqref{eq:discrete_chain_rule} is used, to obtain the stability equality:
\begin{equation} \label{eq:stability_frac_flow_1}
    \begin{aligned}
        &\frac{c_{f}}{2\Delta t}\left(\lVert g^{1/2}\pch \rVert_{L^2(\mc)}^2-\lVert g^{1/2}\pcnp \rVert_{L^2(\mc)}^2+\lVert  g^{1/2}(\pch-\pcnp) \rVert_{L^2(\mc)}^2\right)-([\delta_t\bu_n^h]_{\nu},\pch)_{\mc}\\
        &+\frac{1}{12\mu}\lVert g^{3/2}\tnabla \pch \rVert_{\bltwoc}^2=( g^{1/2}\sfrac_n,g^{1/2}\pch)_{\mc}-\frac{1}{6\mu}\sum_{*\in\pm}\Big(g^{1/2}(\pch-\pnhtr),g^{1/2}\pch\Big)_{\mc}.
    \end{aligned}
\end{equation}
Next, we test the variational inequality \eqref{eq:discrete_mechanics} with $\bv^h=\bm0$, and use the discrete chain rule \eqref{eq:discrete_chain_rule} on the bilinear form $A$, resulting in
\begin{equation} \label{eq:stability_mechanics_1}
    \begin{aligned}
        \frac{1}{2\Delta t}\left(A[\bu_n^h,\bu_n^h]\right.-A[\bu_{n-1}^h,\bu_{n-1}^h]&+A[\left. \bu_n^h-\bu_{n-1}^h, \bu_n^h-\bu_{n-1}^h]\right)+B[\frunh,\frunh]\\
        &\leq \alpha\left(p_n^h,\nabla \cdot \frunh\right)-(\pch,[\frunh]_{\nu})_{\mc}-J[\bu_n^h,\frunh]+\left(\bm{f}_n,\frunh\right).
    \end{aligned}
\end{equation}
Equations \eqref{eq:stability_flow_1}, \eqref{eq:stability_frac_flow_1} and \eqref{eq:stability_mechanics_1} are then added, and we get, noting that the coupling terms between flow and mechanics cancel:
\begin{equation} \label{eq:stability_combined_1}
    \begin{aligned}
        &\frac{1}{2\Delta t}\left(A [\bu_n^h,\bu_n^h]-A[\bu_{n-1}^h,\bu_{n-1}^h]+A[ \bu_n^h-\bu_{n-1}^h, \bu_n^h-\bu_{n-1}^h]\right)+B[\frunh,\frunh]\\
        &+\frac{c_0}{2\Delta t}\left(\lVert p_n^h \rVert_{\ltwo}^2-\lVert p_{n-1}^h \rVert_{\ltwo}^2+\lVert p_n^h-p_{n-1}^h
        \rVert_{\ltwo}^2\right)+\frac{1}{\mu}\lVert \bk^{1/2}\nabla p_n^h \rVert_{\bltwopm}^2\\
        &+\frac{c_{f}}{2\Delta t}\left(\lVert g^{1/2}\pch \rVert_{L^2(\mc)}^2-\lVert g^{1/2}\pcnp \rVert_{L^2(\mc)}^2+\lVert  g^{1/2}(\pch-\pcnp) \rVert_{L^2(\mc)}^2\right)+\frac{1}{12\mu}\lVert g^{3/2}\tnabla \pch \rVert_{\bltwoc}^2\\
        &+\frac{1}{6\mu}\sum_{*\in\pm}\lVert g^{1/2}(\pch-\pnhtr)\rVert_{L^2(\mc)}^2\\
        &\leq \left(\psi_n,p_n^h\right)+( g^{1/2}\sfrac_n,g^{1/2}\pch)_{\mc}+\left(\bm{f}_n,\frunh\right)-J[\bu_n^h,\frunh].
    \end{aligned}
\end{equation}
Henceforth, we will drop the non-negative terms $A[\bu_n^h-\bu_{n-1}^h, \bu_n^h-\bu_{n-1}^h]$, $\lVert p_n^h-p_{n-1}^h\rVert_{\ltwo}^2$, $\lVert g^{1/2}(\pch-\pcnp)\rVert_{L^2(\mc)}^2$ and $\sum_{*\in\pm}\lVert g^{1/2}(\pch-\pnhtr)\rVert_{L^2(\mc)}^2$, as these are not needed in the upcoming estimates. We sum each side of \eqref{eq:stability_combined_1} from $n=1$ to $N$, multiply by $\Delta t$, use the continuity \eqref{eq:continuity_bilinear} and coercivity \eqref{eq:v-ellipticity} of the bilinear forms, and additionally use the upper and lower bounds on the permeability $\bk$ \eqref{eq:spd_perm} and fracture gap $g$ \eqref{eq:g_bound} to get
\begin{equation} \label{eq:stability_almost_final}
    \begin{aligned}
        &\frac{\alpha_A}{2}\lVert\bu_N^h\rVert_{\bhonepm}^2+\frac{c_0}{2}\lVert p_N^h\rVert_{\ltwo}^2+\frac{c_{f}g_{\text{min}}}{2}\lVert \pcnh\rVert_{L^2(\mc)}^2\\
        &+\alpha_B \sum_{n=1}^N \Delta t\lVert\frunh\rVert_{\bhonepm}^2+\frac{\lambda_{\text{min}}}{\mu}\sum_{n=1}^N\Delta t \lVert \nabla p_n^h \rVert_{\bltwopm}^2+\frac{g_{\text{min}}^3}{12\mu}\sum_{n=1}^N\Delta t \lVert \tnabla \pch \rVert_{\bltwoc}^2\\
        &\leq \frac{m_A}{2}\lVert\bu_0^h\rVert_{\bhonepm}^2+\frac{c_0}{2}\lVert p_0^h \rVert_{\ltwo}^2+\frac{c_{f}g_{\text{max}}}{2}\lVert \pcnull^h \rVert_{L^2(\mc)}^2\\
        &+ \Delta t \sum_{n=1}^N\left( \left(\psi_n,p_n^h\right)+( g^{1/2}\sfrac_n,g^{1/2}\pch)_{\mc}+\left( \bm{f}_n,\frunh\right)-J[\bu_n^h,\frunh]\right).
    \end{aligned}
\end{equation}
The terms inside the sum on the right-hand side are subsequently bounded. This follows the same procedure as the stability estimate of \cite{nevland2026numerical}, which is easily adaptable to our fractured setting, so we skip the details of these estimates. We merely remark that most of the estimates consist of simple usage of the Cauchy-Schwarz and Young's inequalities with suitable coefficients, before applying the discrete Grönwall inequality in the end. The only part that requires careful consideration is the estimate on the nonlinear contact term $J[\bu_n^h, \frunh]$, where we obtain, after applying Young's inequality, two terms of the form $\varepsilon\lVert \frunh \rVert_{\bhonepm}^2$ and $(4\varepsilon)^{-1}\lVert \bu_n^h \rVert_{\bhonepm}^2$, where $\varepsilon>0$ is the freely chosen coefficient of Young's inequality. The sizes of these terms cannot be controlled simultaneously, and as a result, we must instead control the time step size when applying Grönwall's inequality, to ensure the strict positivity of the terms. Specifically, we must require
\begin{equation}
    \frac{\alpha_A}{2}-(4\varepsilon)^{-1}\Delta t>0,
\end{equation}
and since we must also require $\varepsilon<\alpha_B$ to absorb the term $\varepsilon\lVert \frunh \rVert_{\bhonepm}^2$ into the left-hand side of \eqref{eq:stability_almost_final}, this leads to the following restriction on the time step size:
\begin{equation}
    \Delta t < 2\alpha_A \alpha_B.
\end{equation}
The final estimate obtained, given this time step restriction, is:
\begin{equation} \label{eq:stability_final}
\begin{aligned}
    &\lVert \bu_N^h\rVert_{\bhonepm}^2+\lVert p_N^h \rVert_{\ltwo}^2+\lVert \pcnh \rVert_{L^2(\mc)}^2\\
    &+\sum_{n=1}^N \Delta t \lVert \delta_t \bu_n^h \rVert_{\bhonepm}^2+ \sum_{n=1}^N \Delta t \lVert \nabla p_n^h \rVert_{\bltwopm}^2+\sum_{n=1}^N \Delta t \lVert \tnabla \pch \rVert_{\bltwoc}^2\\
    &\leq C\Bigg( \lVert \bu_0^h \rVert_{\bhonepm}^2+\lVert p_0^h \rVert_{\ltwo}^2+\lVert \pcnull^h\rVert_{L^2(\mc)}^2\\
    &+\Delta t\sum_{n=1}^N\big(\lVert \psi_n \rVert_{L^2(\Omega)}^2+\lVert \sfrac_n\rVert_{L^2(\mc)}^2+\lVert \bm{f}_n\rVert_{L^2(\Omega)}^2+\lVert g\rVert_{L^2(\mc)}^2\big)\Bigg),
\end{aligned}
\end{equation}
where $C$ depends on the final time $T$, due to the Grönwall inequality, but is independent of $n, h$ and $\Delta t$. Since $\bu_0^h=\bu_0^I$, $p_0^h=p_0^I$ and $\pcnull^h=\pcnull^I$, the first three terms on the right-hand side can be bounded, up to a constant, by $\lVert \bu_0\rVert_{\bhonepm}^2$, $\lVert p_0 \rVert_{\ltwo}^2$ and $\lVert \pcnull\rVert_{L^2(\mc)}^2$. For the final four terms on the right-hand side, since the data is continuous on $[0,T]$, we can bound the Riemann sums, up to a constant, by a continuous-in-time norm. For instance,
\begin{equation}
    \Delta t \sum_{n=1}^N \lVert \psi_n \rVert_{L^2(\Omega)}^2 \leq T\lVert \psi \rVert^2_{L^{\infty}(0,T;L^2(\Omega))}.
\end{equation}
Hence, the entire right-hand side can be bounded by a constant independent of $n,h$ and $\Delta t$.
\end{proof}

\subsection{Error estimate} \label{sec:error_estimates}

In this section, we derive an a priori error estimate between the solution to the fully discrete problem (Q-D) and the (presumed) solution to the continuous problem (Q). This derivation requires relatively strong regularity assumptions on the solution; we will comment on the practical feasibility of these assumptions at the end of this section.

\begin{theorem} \label{thm:error_estimate}
    If $\bu \in H^2(0,T;\bm{H}^{q}(\Omega^{\pm}))$, $p \in H^1(0,T;H^{r}(\Omega^{\pm})) \cap H^2(0,T;L^2(\Omega))$ and $\pc \in H^1(0,T;H^{s}(\mc)) \cap H^2(0,T;L^2(\mc))$, with $q, r$ and $s$ satisfying \eqref{eq:sobolev_exponents}, and the time step size satisfies $\Delta t < 2 \alpha_A \alpha_B$, then the sequence of solutions $(\bu_n^h,p_n^h, \pch)$, $1 \leq n \leq N$, to problem (Q-D), satisfies:
        \begin{equation} \label{eq:asymptotic_estimate}
        \begin{aligned}
            & \max_{1 \leq n \leq N}\lVert \bu_n-\bu_n^h\rVert_{\bhonepm}^2+\max_{1\leq n \leq N}\lVert p_n-p_n^h\rVert_{\ltwo}^2 +\max_{1\leq n \leq N}\lVert \pcn-\pch\rVert_{L^2(\mc)}^2\\
            &+\sum_{n=1}^N \Delta t\lVert \nabla (p_n-p_n^h) \rVert_{\bltwopm}^2+\sum_{n=1}^N \Delta t\lVert \tnabla (\pcn-\pch)\rVert_{\bltwoc}^2\\
            &\leq C((\Delta t)^2 + h^{q-1}+h^{2(q-1)}+h^{2(r-1)}+h^{2(s-1)}),
        \end{aligned}
    \end{equation}
    for a constant $C$ independent of $n, h$ and $\Delta t$.
\end{theorem}
\begin{proof}
The proof is split up into three steps. We first derive an error equality between the continuous and discrete mass balance equations, then we derive an error inequality between the continuous and discrete variational inequalities, and finally the two estimates are combined.

\textbf{Step 1: Flow}

Our strategy for the mass balance equations closely follows that of our previous paper \cite{nevland2026numerical}, which in turn is inspired by the papers of Girault et al. \cite{poroelastic_notmixed,Girault2019}, in which error estimates were derived for a poroelastic model in a fractured medium, but excluding fracture contact mechanics. Starting with the mass balance equation in the matrix, the discrete equation \eqref{eq:discrete_biot} is subtracted from the continuous equation \eqref{eq:weak_biot} evaluated at time $t_n$. Setting $q=q^h$, we obtain:
\begin{equation} \label{eq:biot_estimate_1}
\begin{aligned}
    \left(c_0(\dot{p}_n-\frpnh)+\alpha\nabla \cdot (\dotu_n-\frunh),q^h\right)&+\frac{1}{\mu}\left(\bk\nabla(p_n-p_n^h),\nabla q^h\right)\\
    &=\frac{1}{6\mu}\sum_{*\in\pm}\big(g(\pcn-\pch)-g(\pntr-\pnhtr),\qhtr\big)_{\mc}.
\end{aligned}
\end{equation}
Similarly, we subtract the discrete fracture mass balance equation \eqref{eq:discrete_biot_fracture} from the continuous one \eqref{eq:weak_biot_fracture} evaluated at time $t_n$, and set $\qc=\qch$ to obtain
\begin{equation} \label{eq:biot_frac_estimate_1}
\begin{aligned}
    \left(c_{f}g(\dot{\pi}_n-\frpcnh)-[\dotu_n-\frunh]_{\nu},\qch\right)_{\mc}&+\frac{1}{12\mu}\left(g^{3}\tnabla(\pcn-\pch),\tnabla \qch\right)_{\mc}\\
    &=-\frac{1}{6\mu}\sum_{*\in\pm}\big(g(\pcn-\pch)-g(\pntr-\pnhtr),\qch\big)_{\mc}.
\end{aligned}
\end{equation}
We split the errors into a spatial discretization error and a temporal discretization error, by adding and subtracting $\delta_t p_n$, $\delta_t \bu_n$ and $\delta_t \pcn$ to the appropriate terms. Furthermore, the spatial discretization errors, $\bu_n-\bu_n^h$, $p_n-p_n^h$ and $\pcn-\pch$, are split into interpolation and auxiliary discretization errors, by adding and subtracting $\bu_n^I$, $p_n^I$ and $\pcn^I$. We define
\begin{equation} \label{eq:error_split}
\begin{aligned}
    &\inu=\bu_n-\bu_n^I, \ \ \inp=p_n-p_n^I, \ \ \inpc = \pcn-\pcn^I,\\
    &\enu=\bu_n^I-\bu_n^h, \ \ \enp=p_n^I-p_n^h,
    \ \ \enpc = \pcn^I-\pch.
\end{aligned}
\end{equation}
Using the above notation, \eqref{eq:biot_estimate_1} can be rewritten as
\begin{equation} \label{eq:biot_estimate_2}
\begin{aligned}
    &\left(c_0\delta_t(\inp+\enp)+\alpha\nabla \cdot \delta_t(\inu+\enu),q^h\right)+\frac{1}{\mu}\left(\bk\nabla(\inp+\enp),\nabla q^h\right)\\
    &=\frac{1}{6\mu}\sum_{*\in\pm}\big(g(\inpc+\enpc)-g(\inptr+\enptr),\qhtr\big)_{\mc}+\left(c_0(\delta_t p_n-\dotp_n)+\alpha\nabla \cdot (\frun-\dotu_n),q^h\right),
\end{aligned}
\end{equation}
while \eqref{eq:biot_frac_estimate_1} is rewritten as
\begin{equation} \label{biot_frac_estimate_2}
    \begin{aligned}
        &\left(c_{f}g\delta_t(\inpc+\enpc)-[\delta_t(\inu+\enu)]_{\nu},\qch\right)_{\mc}+\frac{1}{12\mu}\left(g^{3}\tnabla(\inpc+\enpc),\tnabla \qch\right)_{\mc}\\
        &=-\frac{1}{6\mu}\sum_{*\in\pm}\big(g(\inpc+\enpc)-g(\inptr+\enptr),\qch\big)_{\mc}+\left(c_{f}g(\delta_t \pcn-\dot{\pi}_{n})-[\frun-\dotu_n]_{\nu},\qch\right)_{\mc}.
    \end{aligned}
\end{equation}
We choose $q^h=\enp$, $\qch=\enpc$, add equations \eqref{eq:biot_estimate_2} and \eqref{biot_frac_estimate_2}, move all terms involving interpolation errors over on the right-hand side and use the chain rule \eqref{eq:discrete_chain_rule} to obtain:
\begin{equation} \label{eq:error_flow_final}
    \begin{aligned}
        &\frac{c_0}{2\Delta t}\left(\lVert \enp \rVert_{\ltwo)}^2-\lVert e_{n-1}^p\rVert_{\ltwo}^2+\lVert \enp-e_{n-1}^p\rVert_{\ltwo}^2\right)+\alpha\left(\nabla \cdot \delta_t\enu,\enp\right)\\
        &+\frac{c_{f}}{2\Delta t}\left(\lVert g^{1/2}\enpc \rVert_{L^2(\mc)}^2-\lVert g^{1/2}\enppc\rVert_{L^2(\mc)}^2+\lVert g^{1/2}(\enpc-\enppc)\rVert_{L^2(\mc)}^2\right)-([\delta_t\enu]_{\nu},\enpc)_{\mc}\\
        &+\frac{1}{\mu}\lVert \bk^{1/2}\nabla \enp \rVert_{\bltwopm}^2+\frac{1}{12\mu}\lVert g^{3/2}\tnabla \enpc\rVert_{\bltwoc}^2+\frac{1}{6\mu}\sum_{*\in\pm}\lVert g^{1/2}(\enpc-\enptr)\rVert_{L^2(\mc)}^2\\
        &=E_n^*+F_n^*+G_n^*,
    \end{aligned}
\end{equation}
where we have defined
\begin{align}
    E_n^*&=-\big(c_0\delta_t \inp+\alpha\nabla \cdot \delta_t\inu,\enp\big)-\frac{1}{\mu}\big(\bk\nabla\inp,\nabla\enp\big)+\big(c_0(\delta_t p_n-\dotp_n)+\alpha\nabla \cdot (\frun-\dotu_n),\enp\big),\\
        F_n^*&=-\big(c_fg\delta_t\inpc-[\delta_t\inu]_{\nu},\enpc\big)_{\mc}-\frac{1}{12\mu}\big(g^{3}\tnabla \inpc,\tnabla \enpc\big)_{\mc}+\big(c_{f}g(\delta_t \pcn-\dot{\pi}_{n})-[\frun-\dotu_n]_{\nu},\enpc\big)_{\mc},\\
        G_n^*&=\frac{1}{6\mu}\sum_{*\in\pm}\big(g(\inpc-\inptr),\enptr-\enpc\big)_{\mc}.
\end{align}
\textbf{Step 2: Mechanics}

Next, we turn to the variational inequalities. Here, our strategy is again similar to our previous paper \cite{nevland2026numerical}, which in turn is based on ideas from chapter 11 of the textbook of Han and Sofonea \cite{han2002quasistatic}. We add the continuous variational inequality \eqref{eq:weak_mechanics}, evaluated at $t_n$, to the discrete \eqref{eq:discrete_mechanics} variational inequality, and set $\bv=\delta_t \bu_n^h$:
\begin{equation} \label{eq:inequalities_added_1}
    \begin{aligned}
        &A[\bu_n,\frunh-\dotu_n]+A[\bu_n^h,\bv^h-\frunh]+B[\dotu_n,\frunh-\dotu_n]+B[\frunh,\bv^h-\frunh]\\
        &-\alpha\left(p_n,\nabla \cdot (\frunh-\dotu_n)\right)-\alpha\left(p_n^h,\nabla \cdot (\bv^h-\frunh)\right)+\left(\pcn,[\frunh-\dotu_n]_{\nu}\right)_{\mc}+\left(\pch,[\bv^h-\frunh]_{\nu}\right)_{\mc}\\
        &+J[\bu_n,\frunh]-J[\bu_n,\dotu_n]+J[\bu_n^h,\bv^h]-J[\bu_n^h,\frunh]\geq \left(\bm{f}_n,\frunh-\dotu_n+\bv^h-\frunh\right).
    \end{aligned}
\end{equation}
We aim to rewrite \eqref{eq:inequalities_added_1} to obtain terms containing the errors $\bu_n-\bu_n^h, \dotu_n-\frunh$, $p_n-p_n^h$ and $\pcn-\pch$. Doing so will lead to an additional residual-type term, which is typical for variational inequalities \cite{han2002quasistatic,kikuchi_oden}. The terms on the left-hand side of \eqref{eq:inequalities_added_1} are rewritten to obtain:
\begin{equation} \label{eq:inequalities_added_2}
    \begin{aligned}
        &A[\bu_n-\bu_n^h,\dotu_n-\frunh]+B[\dotu_n-\frunh,\dotu_n-\frunh]-\alpha\left(p_n-p_n^h,\nabla \cdot (\dotu_n-\frunh)\right)+\left(\pcn-\pch,[\dotu_n-\frunh]_{\nu}\right)_{\mc}\\
        &\leq A[\bu_n^h-\bu_n,\bv^h-\dotu_n]+B[\frunh-\dotu_n,\bv^h-\dotu_n]-\alpha\left(p_n^h-p_n,\nabla \cdot (\bv^h-\dotu_n)\right)+\left(\pch-\pcn,[\bv^h-\dotu_n]_{\nu}\right)_{\mc}\\
        &+J[\bu_n,\frunh]-J[\bu_n,\bv^h]+J[\bu_n^h,\bv^h]-J[\bu_n^h,\frunh]+R_n(\bv^h),
    \end{aligned}
\end{equation}
where the residual-type term $R_n$ is defined as follows:
\begin{equation} \label{eq:residual}
\begin{aligned}    R_n(\bv^h) = A&[\bu_n,\bv^h-\dotu_n]+B[\dotu_n,\bv^h-\dotu_n]-\alpha\left(p_n,\nabla \cdot (\bv^h-\dotu_n)\right)+\left(\pcn,[\bv^h-\dotu_n]_{\nu}\right)_{\mc}\\
&+J[\bu_n,\bv^h]-J[\bu_n,\dotu_n]-\left(\bm{f}_n,\bv^h-\dotu_n\right).
\end{aligned}
\end{equation}
As for the mass balance equation, we split the error into temporal and spatial discretization errors, and the latter is further split into an interpolation and auxiliary discretization error. More specifically, we write, using the same notation \eqref{eq:error_split} as before:
\begin{align}
    \bu_n-\bu_n^h &= \inu + \enu \\
    \dotu_n-\frunh &= \dotu_n - \frun + \franu + \frenu.
\end{align}
Moreover, we set $\bv^h=\dotu_n^I$, and use the notation $\dot{a}_n^{\bu}=\dotu_n-\dotu_n^I$. This particular test function is chosen because $\dot{a}_n^{\bu}$ can be bounded by the interpolation error estimate \eqref{eq:interpolation_estimate}. After an elementary but lengthy calculation, we can rewrite \eqref{eq:inequalities_added_2} to
\begin{equation} \label{eq:final_mechanics_equation}
    \begin{aligned}
        &\frac{1}{2\Delta t}\left(A[\enu,\enu]-A[e_{n-1}^{\bu},e_{n-1}^{\bu}]+A[\enu-e_{n-1}^{\bu},\enu-e_{n-1}^{\bu}]\right)+B[\frenu,\frenu]\\
        &-\alpha\left(\enp,\nabla \cdot \frenu\right)+(\enpc,[\frenu]_{\nu})_{\mc}\leq A_n^* + B_n^* + C_n^* + D_n^* + J_n^* + R_n(\dotu_n^I),
    \end{aligned}
\end{equation}
where
\begin{equation}
    A^*_n = A[\inu+\enu,\dot{a}_n^{\bu}] -A[\enu,\dotu_n-\frun+\franu]-A[\inu,\dotu_n-\frun+\franu+\frenu],
\end{equation}
\begin{equation}
\begin{aligned}
    B^*_n = B&[\franu+\frenu+\dotu_n-\frun,\dot{a}_n^{\bu}]-B[\dotu_n-\frun,\dotu_n-\frun]-B[\franu,\franu]\\
        &-2B[\dotu_n-\frun,\franu]-2B[\frenu,\franu+\dotu_n-\frun],
\end{aligned}
\end{equation}
\begin{equation}
    C_n^* = \alpha\left(\inp,\nabla \cdot (\dotu_n-\frun+\franu+\frenu)\right)+\alpha\left(\enp,\nabla \cdot (\dotu_n-\frun+\franu)\right)-\alpha\left(\inp+\enp,\dot{a}_n^{\bu}\right),
\end{equation}
\begin{equation}
    D_n^* = -(\inpc,[\dotu_n-\frun+\franu+\frenu]_{\nu})_{\mc}-(\enpc,[\dotu_n-\frun+\franu]_{\nu})_{\mc}+(\inpc+\enpc,[\dot{a}_n^{\bu}]_{\nu})_{\mc},
\end{equation}
\begin{equation}
    J_n^* = J[\bu_n,\frunh]-J[\bu_n,\dotu_n^I]+J[\bu_n^h,\dotu_n^I]-J[\bu_n^h,\frunh],
\end{equation}
and the residual term $R_n$ is defined as before, but with $\bv^h=\dotu_n^I$ substituted. Note that, just like \eqref{eq:error_flow_final}, we have moved all terms involving interpolation operators to the right-hand side of \eqref{eq:final_mechanics_equation}.

\textbf{Step 3: Combining flow and mechanics}

The next step is then to combine \eqref{eq:error_flow_final} and \eqref{eq:final_mechanics_equation}. Noting that the coupling terms $\alpha(\enp,\nabla \cdot \frenu)$ and $(\enpc,[\frenu]_{\nu})_{\mc}$ cancel, we obtain:
\begin{equation} \label{eq:big_estimate}
    \begin{aligned}
        &\frac{1}{2\Delta t}\left(A[\enu,\enu]-A[e_{n-1}^{\bu},e_{n-1}^{\bu}]+A[\enu-e_{n-1}^{\bu},\enu-e_{n-1}^{\bu}]\right)+B[\frenu,\frenu]\\
        &+\frac{c_0}{2\Delta t}\left(\lVert \enp \rVert_{\ltwo}^2-\lVert e_{n-1}^p\rVert_{\ltwo}^2+\lVert \enp-e_{n-1}^p\rVert_{\ltwo}^2\right)+\frac{1}{\mu}\lVert \bk^{1/2}\nabla \enp \rVert_{\bltwopm}^2\\
        &+\frac{c_{f}}{2\Delta t}\left(\lVert g^{1/2}\enpc \rVert_{L^2(\mc)}^2-\lVert g^{1/2}\enppc\rVert_{L^2(\mc)}^2+\lVert g^{1/2}(\enpc-\enppc)\rVert_{L^2(\mc)}^2\right)+\frac{1}{12\mu}\lVert g^{3/2}\tnabla \enpc\rVert_{\bltwoc}^2\\
        &+\frac{1}{6\mu}\sum_{*\in\pm}\lVert g^{1/2}(\enpc-\enptr)\rVert_{L^2(\mc)}^2\\
        &\leq A_n^* + B_n^* + C_n^* + D_n^* + E_n^* + F_n^* + G_n^* + J_n^* + R_n(\dotu_n^I),
    \end{aligned}
\end{equation}
Henceforth, we will drop the non-negative terms $A[e_{n}^{\bu}-e_{n-1}^{\bu},e_{n}^{\bu}-e_{n-1}^{\bu}]$, $\lVert \enp-e_{n-1}^p\rVert_{\ltwo}^2$ and $\lVert g^{1/2}(\enpc-\enppc)\rVert_{L^2(\mc)}^2$ from the left-hand side of \eqref{eq:big_estimate}, as these terms are not needed for the upcoming estimates. Multiplying both sides of \eqref{eq:big_estimate} by $\Delta t$, summing both sides from $n=1$ to $N$, using the coercivity \eqref{eq:v-ellipticity} of the bilinear forms, and using the lower bounds of the permeability $\bk$ \eqref{eq:spd_perm} and fracture gap $g$ \eqref{eq:g_bound}, we get
\begin{equation} \label{eq:big_estimate_summed}
    \begin{aligned}
        &\frac{\alpha_A}{2}\lVert e_N^{\bu}\rVert_{\bhonepm}^2+\frac{c_0}{2}\lVert e_N^p\rVert_{\ltwo}^2+\frac{c_{f}g_{\text{min}}}{2}\lVert e_{N}^{\pi}\rVert_{L^2(\mc)}^2+\alpha_B\sum_{n=1}^N\Delta t\left\lVert \frenu \right\rVert_{\bhonepm}^2\\
        &+\frac{\lambda_{\text{min}}}{\mu}\sum_{n=1}^N\Delta t\lVert \nabla \enp \rVert_{\bltwopm}^2+\frac{g_{\text{min}}^3}{12\mu}\sum_{n=1}^N \Delta t\lVert \tnabla \enpc\rVert_{\bltwoc}^2+\frac{g_{\text{min}}}{6\mu}\sum_{n=1}^N\Delta t\sum_{*\in\pm}\lVert \enpc-\enptr\rVert_{L^2(\mc)}^2\\
    &\leq \Delta t \sum_{n=1}^N \lvert A_n^* \rvert + \lvert B_n^* \rvert + \lvert C_n^* \rvert + \lvert D_n^* \rvert+  \lvert E_n^* \rvert + \lvert F_n^* \rvert  + \lvert G_n^* \rvert +\lvert J_n^* \rvert + \lvert R_n(\dotu_n^I) \rvert.
    \end{aligned}
\end{equation}
Note that due to our choice of initial conditions, the terms $A[e_0^{\bu},e_0^{\bu}]$, $\lVert e_0^p\rVert_{\ltwo}^2$ and $\lVert e_0^{\pi}\rVert_{L^2(\mc)}^2$ vanish. It remains to bound the terms on the right-hand side of \eqref{eq:big_estimate_summed}. Since this involves using quite simple but repetitive techniques, we move these estimates to Appendix \ref{app:error}. After bounding all terms on the right-hand side and applying the discrete Grönwall inequality (which requires the restriction $\Delta t < \alpha_A\alpha_B$ on the time step, for the same reason as in the stability estimate), we obtain the following error estimate:
\begin{equation} \label{eq:final_estimate}
    \begin{aligned}
        &\lVert e_N^{\bu}\rVert_{\bhonepm}^2+\lVert e_N^p\rVert_{\ltwo}^2+\lVert e_{N}^{\pi}\rVert_{L^2(\mc)}^2+\sum_{n=1}^N\Delta t\left\lVert \frenu \right\rVert_{\bhonepm}^2+\sum_{n=1}^N\Delta t\lVert \nabla \enp \rVert_{\bltwopm}^2\\
    &+\sum_{n=1}^N\Delta t\lVert \tnabla \enpc \rVert_{\bltwoc}^2 \leq C \Bigg( \lVert \dotu-\dotu^I\rVert_{L^2(0,T;\bhonepm)}^2 + \lVert \dotp-\dotp^I\rVert_{L^2(0,T;\ltwo)}^2\\
    &+\lVert \dotpc-\dotpc^I\rVert_{L^2(0,T;L^2(\mc))}^2+ (\Delta t)^2\left(\lVert \ddot{\bu}\rVert_{L^2(0,T;\bhonepm)}^2+\lVert \ddot{p}\rVert_{L^2(0,T;\ltwo)}^2+\lVert \ddotpc\rVert_{L^2(0,T;L^2(\mc))}^2\right)\\
    &+\Delta t \sum_{n=1}^N \Big(\lVert \dot{a}_n^{\bu}\rVert_{\bhonepm}+\lVert \inu \rVert_{\bhonepm}^2 + \lVert \dot{a}_n^{\bu} \rVert_{\bhonepm}^2+\lVert \inp \rVert_{\honepm}^2 + \lVert \inpc \rVert_{H^1(\mc)}^2\Bigg),
    \end{aligned}
\end{equation}
where $C$ depends on the final time $T$, but is independent of $n, h$ and $\Delta t$.
Like in the stability estimate, because $\bu,\dotu,p$ and $\pc$ are continuous on $[0,T]$, due to our assumptions on their temporal regularity, the Riemann sums on the right-hand side of the inequality can be bounded, up to a constant, by a continuous-in-time norm. For instance:
\begin{equation}
    \Delta t \sum_{n=1}^N \lVert \dot{a}_n^{\bu}\rVert_{\bhonepm} \leq T \lVert \dotu-\dotu^I \rVert_{L^{\infty}(0,T;\bhonepm)}.
\end{equation}
Finally, by employing the interpolation error estimates \eqref{eq:interpolation_estimate} on the relevant terms on the right-hand side of \eqref{eq:final_estimate}, we obtain the desired asymptotic error estimate \eqref{eq:asymptotic_estimate}.
\end{proof}

\begin{remark} \label{rem:regularity}
Let us comment on the practical feasibility of the regularity assumptions of Theorem \ref{thm:error_estimate}, specifically in the context of using linear elements, which will be the case for the numerical experiment of the next Section. For the Biot equations in non-fractured media, arbitrary spatial regularity can be obtained, given sufficiently smooth data, from standard elliptic and parabolic regularity theory \cite{showalter2000diffusion}. The situation is more delicate with contact mechanics involved, as the contact conditions tend to reduce spatial regularity, and the optimal assumption (for linear elements) of $H^2$-regularity in space might not hold. It is common to assume the spatial regularity for contact problems to be $H^s, s>3/2$ \cite{banz2025contact,chouly2023finite}. Even more important to consider is the presence of fracture tips. It is well known that, for elliptic problems, the presence of corners/tips in the domain lead to singular expansions near the corner/tip and reduced regularity of the solution \cite{costabel2002crack,grisvard}. Specifically, it is shown in \cite{costabel2002crack} that only a regularity of $H^s, 1 < s <3/2 $ can be expected due to the fracture tips. However, all of these cases are included in Theorem \ref{thm:error_estimate}, which makes the very flexible assumption of the spatial regularity being $H^{s}, s\in(1,2]$, for all variables, which can be expected to be fulfilled even with singularities near the fracture tips.

Regarding the temporal regularity, we note that in chapter 11 of \cite{han2002quasistatic}, a viscoelastic contact problem with normal compliance and Coulomb friction is considered, and it is proven that under standard continuity and coercivity assumptions on the variational inequality, the solution is $H^2$ in time if the data is $H^1$ in time. Again for the Biot equations in non-fractured media, arbitrary regularity in time can be obtained, given sufficiently smooth data, for $t>0$ \cite{showalter2000diffusion}. One must be more careful with $t=0$, as one can not expect arbitrary regularity for all $t \geq 0$ if the initial condition is not compatible with parameters and data. However, if the solution is considered only on a strictly positive time interval, then our assumption in Theorem \ref{thm:error_estimate} of $H^2$-regularity in time might be reasonable.
\end{remark}

\section{Numerical experiment} \label{sec:num_exp}

In this section, we present a numerical experiment in which numerical convergence orders are computed and compared to the theoretical error estimate \eqref{eq:asymptotic_estimate}. The numerical experiment has been implemented in multiphenicsx \cite{multiphenicsx}, an extension of the finite element software FEniCSx \cite{fenics4,fenics1,fenics3,fenics2} that facilitates the restriction of variables to specific subdomains/boundaries. The code for the numerical experiment was developed with partial assistance of the generative AI tool ChatGPT. All suggestions by AI were carefully reviewed and verified by the corresponding author.

Before describing the simulation setup, some technical details of the implementation will first be outlined. The continuity conditions \eqref{eq:interface_cond_flow}, \eqref{eq:interface_cond_mech} across $\Gamma\setminus\mc$ are enforced by the use of Lagrange multipliers. A penalty method is employed for the mechanics subproblem, leading to an algorithm similar to the one in Section 3.2 of \cite{SIMO199297}. The penalty method turns the discrete variational inequality \eqref{eq:discrete_mechanics} into an equality, and the resulting nonlinear equation system is solved by Newton's method together with a backtracking line search, using absolute and relative tolerances of $10^{-8}$. The nonlinear systems quickly get ill-conditioned as the penalty parameter increases, eventually leading to convergence issues. The penalty parameter has accordingly been chosen as large as possible without causing convergence issues; the precise value varies a bit between the different experiments. The linear systems are solved by direct LU factorization.

We perform two separate convergence analyses, one with respect to the spatial discretization parameter $h$
and one with respect to the temporal discretization parameter $\Delta t$. As no analytical solution is available, we use a numerical solution produced with a finer discretization as our reference solution. For the spatial
convergence analysis, this corresponds to a numerical solution produced on a sufficiently fine spatial grid,
while for the temporal convergence analysis, our reference solution is produced by taking sufficiently small
time steps. Denote the reference solutions at time step $n$ by $\bu_n^{\text{ref}}, p_n^{\text{ref}}, \pcn^{\text{ref}}$. In accordance with the theoretical error estimate \eqref{eq:asymptotic_estimate}, the errors of displacement, matrix pressure and fracture pressure are computed as follows:
\begin{subequations} \label{eq:h1_errors}
\begin{align}
    &\text{Error} \ \bu = \max_{1 \leq n \leq N}\lVert \bu_n^{\text{ref}}-\bu_n^h\rVert_{\bhonepm},\\
    &\text{Error} \ p = \max_{1\leq n \leq N}\lVert p_n^{\text{ref}}-p_n^h\rVert_{\ltwo} + \left(\Delta t\sum_{n=1}^N \lVert \nabla (p_n^{\text{ref}}-p_n^h) \rVert_{\bltwopm}^2\right)^{1/2},\\
    &\text{Error} \ \pc = \max_{1\leq n \leq N}\lVert \pcn^{\text{ref}}-\pch\rVert_{L^2(\mc)} + \left(\Delta t\sum_{n=1}^N \lVert \nabla_{\tau} (\pcn^{\text{ref}}-\pch) \rVert_{\bltwoc}^2\right)^{1/2},
\end{align}
\end{subequations}
We proceed to describe our simulation setup, which is designed to produce a complex solution with tightly coupled processes. Our fractured domain is the unit square, $\Omega=(0,1) \times (0,1)$, and we consider a horizontal fracture $\mc=(0.3,0.7) \times \{0.5\}$. The interface separating $\Omega^+$ and $\Omega^-$ is the horizontal line containing the fracture, $\Gamma=[0,1]\times\{0.5\}$, and we let $\Omega^+=(0,1)\times(0.5,1)$ and $\Omega^-=(0,1)\times(0,0.5)$. The external boundary $\partial \Omega$ is divided as follows: $\partial \Omega_d=\partial \Omega_p=[0,1]\times \{0\}$, $\partial\Omega_f=\partial\Omega\setminus\partial\Omega_p$ and $\partial\Omega_n=\partial\Omega \setminus \partial \Omega_d$. Initial values of pressure and displacement are set to zero, and the reference fracture gap is set to a constant value of $g=10^{-4}$ in $\mc$. 
The boundary conditions are homogeneous except from the top boundary, where we assign $\bm{\sigma}^{\text{por}}\cdot\bm{\nu}=(0.001,-0.1)$. The inhomogeneous Neumann condition has two purposes; first, the compressive traction is used to (partially) mechanically close the fracture and induce nontrivial contact behavior, and second, the shear traction is intended to induce stick-slip behavior of the fracture. However, the shear traction is deliberately chosen to not be sufficiently large to overcome the friction bound. Instead, to induce slip on parts of the fracture, we set the source term in the fracture to be $\sfrac=1.5\cdot10^2$, thereby increasing the fracture pressure and  reducing the effective normal contact traction. Unitary values are chosen for all material parameters, with the exception of the viscosity, which is set to $\mu = 10^{-4}$; the reason is that unless the viscosity is sufficiently low, the added pressure that results from compressing the fracture gap prevents the fracture from closing.

\begin{figure}[H]
    \centering
    \captionsetup[subfigure]{labelformat=empty}
    \begin{subfigure}{0.44\textwidth} 
    \includegraphics[width=\textwidth]{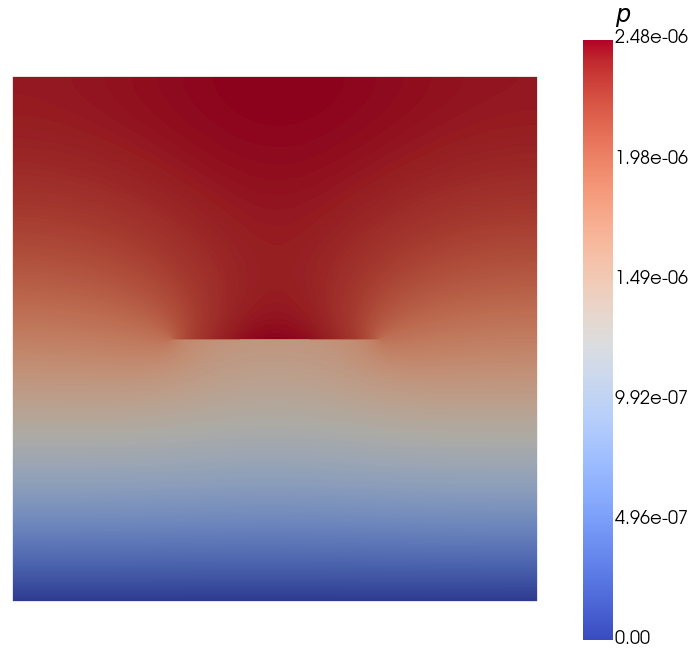}
    \caption{}
    \label{subfig:idk1}
    \end{subfigure}
    \begin{subfigure}{0.55\textwidth} 
    \includegraphics[width=\textwidth]{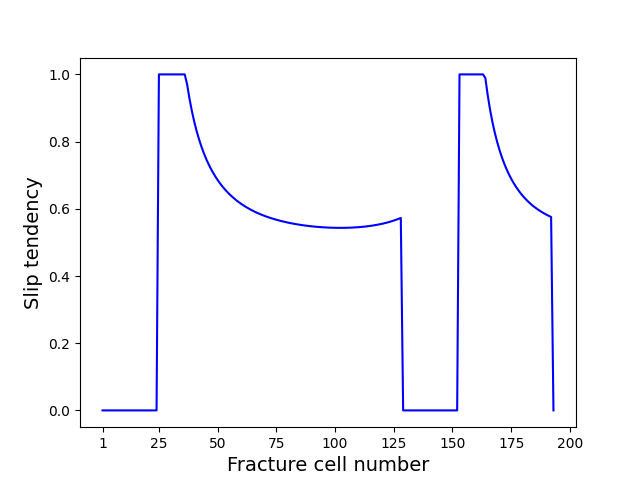}
    \caption{}
    \label{subfig:idk2}
    \end{subfigure}
    \caption{Visualization of the reference solution of the spatial convergence analysis, showing the matrix pressure field (left) and the slip tendency (right).}    \label{fig:visualization}
\end{figure}

A visualization of the simulation is shown in Figure \ref{fig:visualization}, showing the reference solution corresponding to the spatial convergence analysis detailed below. The left figure shows the matrix pressure field (where we arbitrarily choose to plot the pressure on the "+" side at the fracture nodes), and a pressure discontinuity across the fracture is clearly seen. The right figure shows the so-called slip tendency of each fracture cell, which is defined as $\lvert\bm{t}_{\tau}\rvert/b$, with $b=\mu_{\text{fr}}c_{\nu}([\bm{u}]_{\nu}-g)$ denoting the friction bound. According to the Coulomb friction law \eqref{eq:friction_law}, it is equal to one if the fracture cell slips, and between zero and one if it sticks; we define by convention a value of zero to mean an open fracture cell. We observe a complex combination of both open, sticking and slipping cells.

\begin{figure}[H]
    \centering
    \scalebox{0.65}
    {\includegraphics[width=\textwidth]{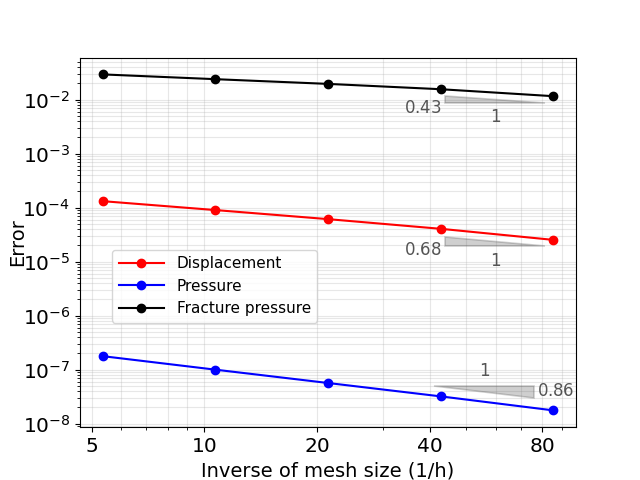}}
    \caption{Errors for various values of the spatial discretization parameter $h$.}
    \label{fig:h1_errors}
\end{figure}

Proceeding with the spatial convergence analysis, we take a single time step with $T=0.01$, and compute the errors on a sequence of nested, refined grids. The reference solution is subsequently produced by performing two additional refinements. The errors \eqref{eq:h1_errors} are computed by interpolating the coarse solutions into the space of the reference solution. The penalty parameter is chosen to be $10^2$ across all grid sizes, and the same value is also used for the constant $c_{\nu}$ of the normal compliance law \eqref{eq:normal_compliance}. The results are shown in Table \ref{tab:errors_space} and Figure \ref{fig:h1_errors}. None of the errors display optimal first order convergence, with particularly the fracture pressure rate being significantly reduced. This is as expected for the displacement, as the theoretical error estimate \eqref{eq:asymptotic_estimate} only implies a convergence order of 0.5 or higher, given that $\bu \in H^2(\Omega^{\pm})$ (the order may even be worse if $\bu$ is only in $H^{2-\varepsilon}$, for $\varepsilon\in (0,1)$). However, the theoretical estimate also implies that the pressure rates should be at least first order if $p\in H^2(\Omega^{\pm})$ and $\pc\in H^2(\mc)$, indicating that this regularity is not achieved in these examples.
\begin{table}[ht!]
\centering
\begin{tabular}{ | l | l | l | l | l | l | l | }
\hline
$1/h$ & Error $p$ & Rate $p$ & Error $\bu$ & Rate $\bu$ & Error $\pc$ & Rate $\pc$ \\
\hline
5.3485 & 1.7814e-07 & - & 1.3180e-04 & - & 2.9258e-02 & - \\
\hline
10.697 & 1.0015e-07 & 0.83089 & 9.0183e-05 & 0.54746 & 2.3875e-02 & 0.29330 \\
\hline
21.394 & 5.6706e-08 & 0.82053 & 6.1231e-05 & 0.55860 & 1.9525e-02 & 0.29018 \\
\hline
42.788 & 3.2055e-08 & 0.82296 & 4.0543e-05 & 0.59481 & 1.5555e-02 & 0.32794 \\
\hline
85.575 & 1.7632e-08 & 0.86232 & 2.5274e-05 & 0.68181 & 1.1558e-02 & 0.42853 \\
\hline
\end{tabular}
\caption{Errors, computed by the formulas \eqref{eq:h1_errors}, and spatial convergence rates.}
\label{tab:errors_space}
\end{table}

\begin{figure}[H]
    \centering
    \scalebox{0.65}
    {\includegraphics[width=\textwidth]{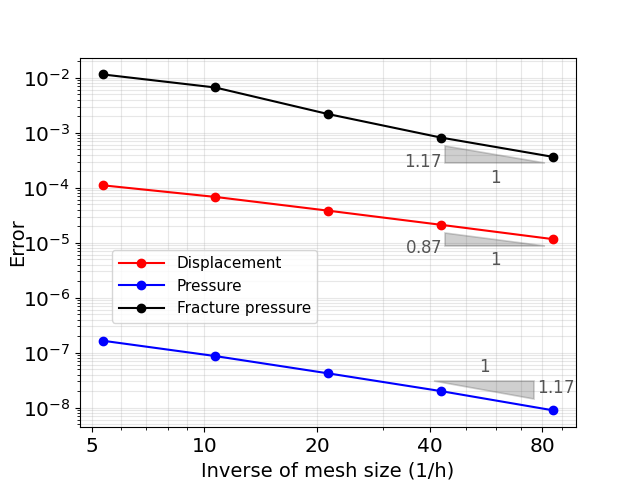}}
    \caption{Errors for various values of the spatial discretization parameter $h$, in the case when a small radius around each fracture tip is excluded in the error norms \eqref{eq:h1_errors}.}
    \label{fig:h1_errors_not_tips}
\end{figure}

As previously mentioned, the solution can be expected to have reduced regularity near the fracture tips, which is likely to lead to reduced convergence rates. It is therefore interesting to recompute the errors of the previous simulation, where balls of a small radius $\delta>0$ around the fracture tips are excluded from the error norms \eqref{eq:h1_errors}. The result is shown in Figure \ref{fig:h1_errors_not_tips}, using $\delta=0.05$. First order convergence is now achieved for the pressure variables, while the order for the displacement is still slightly below one. This indicates that the suboptimal pressure rates observed in Figure \ref{fig:h1_errors} is due to the reduced regularity at the fracture tips.

\begin{table}[ht!]
\centering
\begin{tabular}{ | l | l | l | l | l | l | l | }
\hline
$1/h$ & Error $p$ & Rate $p$ & Error $\bu$ & Rate $\bu$ & Error $\pc$ & Rate $\pc$ \\
\hline
5.3485 & 1.6507e-07 & - & 1.1123e-04 & - & 1.1597e-02 & - \\
\hline
10.697 & 8.6638e-08 & 0.92999 & 6.8283e-05 & 0.70400 & 6.6997e-03 & 0.79165 \\
\hline
21.394 & 4.2023e-08 & 1.0438 & 3.8411e-05 & 0.83002 & 2.2028e-03 & 1.6047 \\
\hline
42.788 & 1.9982e-08 & 1.0725 & 2.1252e-05 & 0.85389 & 8.2289e-04 & 1.4206 \\
\hline
85.575 & 8.8888e-09 & 1.1686 & 1.1608e-05 & 0.87250 & 3.6465e-04 & 1.1742 \\
\hline
\end{tabular}
\caption{Errors and spatial convergence rates in the case when a small radius around each fracture tip is excluded in the error norms \eqref{eq:h1_errors}.}
\label{tab:errors_space_no_tips}
\end{table}

\begin{figure}[H]
    \centering
    \scalebox{0.65}
    {\includegraphics[width=\textwidth]{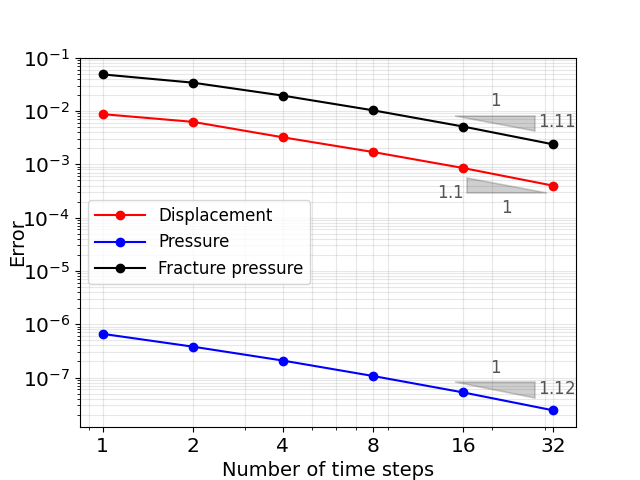}}
    \caption{Errors for various values of the time step $\Delta t$.}
    \label{fig:h1_errors_time}
\end{figure}

Next, we perform a temporal convergence analysis, keeping the spatial mesh constant with $1 /h\approx 48$. In this case, we run the simulation for a longer period of time, until $T=1$; because of the increased ill-conditioning resulting from taking larger time steps, we reduce the penalty parameter to $10^1$. The number of time steps $N$ are varied as $N=2^i, 0\leq1\leq5$. The reference solution is computed with $N=2^8$, and the errors \eqref{eq:h1_errors} are computed by sampling the reference solution at the discrete time points it shares with the coarser solution. Recalling our Remark \ref{rem:regularity}, we consider only time points satisfying $t_n>0.1$ for our error measurements, since we have not designed the initial condition to be compatible with our given parameters and data. As a result, we can expect the solution to initially have a poor temporal regularity, which could adversely affect the convergence rate. The results are presented in Table \ref{tab:errors_time} and Figure \ref{fig:h1_errors_time}. We observe convergence of at least order one in all cases.

\begin{table}[ht!]
\centering
\begin{tabular}{ | l | l | l | l | l | l | l | }
\hline
$T/\Delta t$ & Error $p$ & Rate $p$ & Error $\bu$ & Rate $\bu$ & Error $\pc$ & Rate $\pc$ \\
\hline
1 & 6.6018e-07 & - & 8.7435e-03 & - & 4.8658e-02 & - \\
\hline
2 & 3.7985e-07 & 0.79744 & 6.2677e-03 & 0.48029 & 3.4074e-02 & 0.51399 \\
\hline
4 & 2.0809e-07 & 0.86821 & 3.2112e-03 & 0.96482 & 1.9461e-02 & 0.80807 \\
\hline
8 & 1.0685e-07 & 0.96158 & 1.7014e-03 & 0.91636 & 1.0311e-02 & 0.91637 \\
\hline
16 & 5.2490e-08 & 1.0255 & 8.5158e-04 & 0.99853 & 5.1115e-03 & 1.0124 \\
\hline
32 & 2.4175e-08 & 1.1185 & 3.9671e-04 & 1.1021 & 2.3652e-03 & 1.1118 \\
\hline
\end{tabular}
\caption{Errors, computed by the formulas \eqref{eq:h1_errors}, and temporal convergence rates.}
\label{tab:errors_time}
\end{table}

\section{Conclusion} \label{sec:conclusion}

We have performed a space-time numerical analysis of a mixed-dimensional model for a fractured poroviscoelastic medium including fracture contact mechanics. The resulting variational problem consists of a coupled system of two linear equations (fluid mass balance in the matrix and fracture) and a nonlinear variational inequality (momentum balance of the solid subject to frictional contact conditions). The problem was discretized in space by conformal finite elements, and in time by the implicit Euler method. Our analysis consisted of
proving existence and uniqueness of the discrete solution, deriving stability estimates of the discrete solution,
and deriving an a priori error estimate between the continuous and discrete solutions. The error estimate was derived in a combined space-time energy norm for the sum of the displacement and pressure variables, and the convergence order was optimal in time, while in space the leading order term was reduced by a square root. This reduction of the spatial convergence order by a square root is typical of contact problems.

A numerical experiment was conducted, in which convergence rates were computed and compared to the theoretical error estimate, using linear elements. The convergence analysis was done separately with respect to the spatial and temporal discretization parameters. For the matrix and fracture pressures, the low regularity of the solution around the fracture tips caused a reduction in the spatial convergence rates; this was confirmed by the recovery of optimal first-order rates if a small radius around each fracture tip was excluded in the error computation. For the displacement, the spatial rate did not reach first order, even with the fracture tips removed. For the temporal convergence analysis, optimal first-order rates were obtained for all variables.

The multiphysics problem considered in this work can be extended in several directions. The permeability, both in the matrix and the fracture, has been taken as constant in our linearized model. However, more advanced models, such as the bulk permeability being a nonlinear function of $\nabla \cdot \bu$ \cite{hommel2018porosity}, and the fracture permeability being a cubic function of the width, will in many cases better represent the physics. Extension to these cases remain open. Moreover, in the numerical analysis setting, the optimal $L^2$ error estimates for contact mechanics problems is a problem of ongoing investigation, as the usual duality methods do not work; some early results in this direction can be found in \cite{christof2020finite,steinbach2016trace}. In our setting, this will require a careful consideration of the regularity of the free boundary problem between the open and closed regions of the fracture.

\section*{Acknowledgments}

This project has received funding from the European Research Council (ERC) under the European Union’s
Horizon 2020 research and innovation program (grant agreement No. 101002507).

\appendix

\section{Contraction estimate of the fixed-stress splitting scheme} \label{app:contraction}

In this appendix, we prove Lemma \ref{lemma:fixed_stress_limits}, establishing the existence of limit functions $p_n^h, \pch$ and $\bu_n^h$ of the fixed-stress iterations \eqref{eq:fixed_stress_biot}, \eqref{eq:fixed_stress_biot_fracture} and \eqref{eq:fixed_stress_mechanics} as $k\to \infty$, and proving that these limit functions solve problem (Q-D-S). Our convergence proof of the fixed-stress iterations combines techniques from two papers. The first is the paper of Girault et al. \cite{girault_splitting}, where the fixed-stress scheme was analyzed for a mixed-dimensional model of coupled flow and mechanics, but excluding contact mechanics; the analysis of this paper was conducted on the continuous model. The second is the paper of Almani and Kumar \cite{fixed_stress_contact}, where the fixed-stress scheme was analyzed for the discretized Biot equations (in a non-fractured medium) subject to frictionless contact constraints.

We start by deriving an estimate on the differences in iterations of the splitting scheme, which will lead to a contraction result. We denote differences in iterations of the splitting scheme by $\telta$, to distinguish it from the notation $\delta_t$ used to denote differences between time steps. Hence, for some variable $v^{h}_{n,k}$, depending on the time step $t_n$, the mesh $\mathcal{T}^h$ and the splitting scheme iteration $k$, we define:
\begin{equation} \label{eq:difference_iterates_splitting}
    \telta v^{h}_{n,k}=v^{h}_{n,k}-v^{h}_{n,k-1}.
\end{equation}
Following standard nomenclature, we define the volumetric mean stress
\begin{equation} \label{eq:mean_stress}
    \sigma_{n,k}^h=\lambda\nabla \cdot \upfs-\alpha \ppfs,
\end{equation}
Similarly, we define a "fracture stress" $\xi$, following \cite{girault_splitting}, as follows:
\begin{equation} \label{eq:mean_stress_fracture}
    \xi_{n,k}^h=\frac{\gamma_c}{\chi} \pcpfs-\frac{1}{\chi}\wpfs ,
\end{equation}
where $\chi>0$ is a constant that will be chosen during the upcoming contraction proof. For notational convenience, we also define the following constant:
\begin{equation} \label{eq:splitting_constants}
    \beta=\frac{c_0}{\alpha^2}+\frac{1}{\lambda}.
\end{equation}
\begin{lemma} \label{lemma:fixed_stress_contraction}
    Given a sufficiently small time step size $\Delta t$, the iterations of the fixed-stress splitting scheme \eqref{eq:fixed_stress_biot}, \eqref{eq:fixed_stress_mechanics} satisfy the following estimate:
\begin{equation} \label{eq:final_inequality}
\begin{aligned}
    &\left\lVert \telta \sigma_{n,k+1}^h\right\rVert_{\ltwopm}^2+\left\lVert \telta \xi_{n,k+1}^h\right\rVert_{L^2(\mathcal{C})}^2+\left(4GD\lambda-\frac{1}{\chi^2}+\frac{2\gamma D\lambda}{\Delta t}\right)\left\lVert \telta \wfs \right\rVert_{L^2(\mathcal{C})}^2\\
    &+\left(\lambda^2+\frac{2\gamma\lambda}{\Delta t}\right)\left\lVert \nabla \cdot \telta \ufs\right\rVert_{\ltwopm}^2+\frac{2\lambda_{\text{min}}\Delta t}{\beta\mu}\left\lVert \nabla\telta \pfs\right\rVert_{\bltwopm}^2+\frac{g_{\text{min}}^3\Delta t}{6\beta\mu}\left\lVert \tnabla\telta \pcfs\right\rVert_{\bltwoc}^2\\
    &+\frac{g_{\text{min}}}{3\beta\mu}\sum_{*\in\pm}\lVert \telta \pcfs-\telta \pfstr\rVert_{L^2(\mc)}^2\leq C\left(\left\lVert \telta \sigma_{n,k}^h\right\rVert_{\ltwopm}^2+\left\lVert \telta \xi^{h}_{n,k} \right\rVert_{L^2(\mc)}^2\right),
\end{aligned}
\end{equation}
for a constant $C<1$.
\end{lemma}

\begin{proof}
We start by computing the difference between the matrix flow equations \eqref{eq:fixed_stress_biot} at iterations $k$ and $k+1$, testing with $q^h=\telta \pfs$ in both cases. Additionally, we use the boundedness of the permeability \eqref{eq:spd_perm}. The result is, after multiplying both sides by $\Delta t$:
\begin{equation} \label{fs_matrix_flow1}
    \begin{aligned}
       &\beta\lVert \alpha\telta \pfs\rVert^2_{\ltwo}+\frac{\lambda_{\text{min}}\Delta t}{\mu}\lVert \nabla\telta \pfs\rVert_{\bltwopm}^2\\
       & = \frac{1}{6\mu}\sum_{*\in\pm}\big(g(\pcfs-\pfstr),\pfstr\big)-(\frac{\alpha}{\lambda}\telta \sigma^{h}_{n,k},\telta \pfs)\\
       &\leq \frac{1}{6\mu}\sum_{*\in\pm}\big(g(\pcfs-\pfstr),\pfstr\big)+\frac{\varepsilon}{2}\lVert \alpha \telta \pfs\rVert_{\ltwo}^2+\frac{1}{2\varepsilon \lambda^2}\rVert \telta \sigma^{h}_{n,k}\rVert_{\ltwopm}^2
    \end{aligned}
\end{equation}
where the Cauchy-Schwarz and Young's inequalities have been used to obtain the third line. Similarly, we compute the difference between the fracture flow equation \eqref{eq:fixed_stress_biot_fracture} at iterations $k$ and $k+1$, test with $\qch = \telta \pcfs$ in both cases, use the boundedness \eqref{eq:g_bound} of the fracture gap, and multiply by $\Delta t$ to get:
\begin{equation} \label{eq:fs_fracture_flow1}
    \begin{aligned}
        \beta_c\lVert \telta \pcfs\lVert_{L^2(\mc)}^2&+\frac{g_{\text{min}}^3\Delta t}{12\mu}\lVert \tnabla\telta \pcfs\rVert_{\bltwoc}^2\\
        &\leq-\frac{1}{6\mu}\sum_{*\in\pm}\big(g(\pcfs-\pfstr),\pcfs\big)+(\gamma_c\telta \pcpfs+[\telta \upfs]_{\nu},\telta \pcfs)_{\mc}.
    \end{aligned}
\end{equation}
Where we have defined $\beta_c=c_fg_{\text{min}}+\gamma_c$. Adding \eqref{fs_matrix_flow1} and \eqref{eq:fs_fracture_flow1}, choosing $\varepsilon=\beta$ and multiplying by $2/\beta$, we get
\begin{equation} \label{eq:fixed_stress_estimate_flow}
    \begin{aligned}
        &\lVert \alpha\telta \pfs\rVert^2_{\ltwo}+\frac{2\lambda_{\text{min}}\Delta t}{\beta\mu}\lVert \nabla\telta \pfs\rVert_{\bltwopm}^2+\frac{2\beta_c}{\beta}\lVert \telta \pcfs\lVert_{L^2(\mc)}^2+\frac{g_{\text{min}}^3\Delta t}{6\beta\mu}\lVert \tnabla\telta \pcfs\rVert_{\bltwoc}^2\\
        &+\frac{g_{\text{min}}}{3\beta\mu}\sum_{*\in\pm}\lVert \telta \pcfs-\telta \pfstr\rVert_{L^2(\mc)}^2\leq \frac{1}{(\beta\lambda)^2}\lVert \telta \sigma_{n,k}^h \rVert_{\ltwopm}^2+\frac{2}{\beta}(\gamma_c\telta \pcpfs+[\telta \upfs]_{\nu},\telta \pcfs)_{\mc}.
    \end{aligned}
\end{equation}
Next, we consider the variational inequality \eqref{eq:fixed_stress_mechanics}. In this case, we test \eqref{eq:fixed_stress_mechanics} for iteration $k$ with $\bm{v}^h=\delta_t\ufs$, and for iteration $k+1$ with $\bm{v}^h=\delta_t\upfs$, and add the two inequalities. Since the first argument of the $J$-functionals are fixed, this choice of test functions causes the $J$-functionals to cancel when the inequalities are added. The resulting estimate becomes, after substituting the definitions \eqref{eq:bilinear_form_1} and \eqref{eq:bilinear_form_2} of the bilinear forms:
\begin{equation} \label{eq:fixed_stress_estimate_mechanics}
\begin{aligned}
    (2G+\frac{\gamma}{\Delta t})\lVert \bm{\varepsilon}(\telta \ufs)\rVert_{\bltwopm}^2&+(\lambda+\frac{\gamma}{\Delta t})\lVert \nabla \cdot \telta \ufs\rVert_{\ltwopm}^2\\
    &-\alpha(\telta \pfs,\nabla \cdot \telta \ufs)+(\telta \pcfs,[\telta \ufs]_{\nu})_{\mc}\leq 0.
\end{aligned}
\end{equation}
The final two terms on the left-hand side of \eqref{eq:fixed_stress_estimate_mechanics} can be absorbed into norms by completing the squares. However, in order to do this, we need a term involving the squared norm of $[\telta \ufs]_{\nu}$. This can be done by deriving an inequality between $\bm{\varepsilon}(\bm{u})$ and $[\bm{u}]_{\nu}$. As shown in \cite{girault_splitting}, the following inequality holds:
\begin{equation} \label{eq:eps_jump_inequality}
    \lVert \bm{\varepsilon}(\bm{u})\rVert_{\bltwopm}^2 \geq D\lVert [\bm{u}]_{\nu} \rVert_{L^2(\mathcal{C})}^2.
\end{equation}
The proof of this inequality is a straightforward application of the trace, Poincaré and Korn inequalities, hence the constant $D$ depends on the constants from these inequalities. We substitute this inequality into \eqref{eq:fixed_stress_estimate_mechanics}, and multiply \eqref{eq:fixed_stress_estimate_mechanics} by $2\lambda$ in order to complete the square on the second-to-last term on the left-hand side. In addition, we shall make the substitution $\telta \wfs=-[\telta \ufs]_{\nu}$, in order to flip the sign on the last term on the left-hand side of \eqref{eq:fixed_stress_estimate_mechanics}. The result is the alternative estimate:
\begin{equation} \label{eq:fixed_stress_estimate_mechanics_2}
\begin{aligned}
    &(4GD\lambda+\frac{2\gamma D\lambda}{\Delta t})\lVert \telta \wfs\rVert_{L^2(\mathcal{C})}^2+(2\lambda^2+\frac{2\gamma\lambda}{\Delta t})\lVert \nabla \cdot \telta \ufs\rVert_{\ltwopm}^2\\
    &-2\alpha\lambda(\telta \pfs,\nabla \cdot \telta \ufs)-2\lambda(\telta \pcfs,\telta \wfs)_{\mc}\leq 0.
\end{aligned}
\end{equation}
Adding \eqref{eq:fixed_stress_estimate_flow} and \eqref{eq:fixed_stress_estimate_mechanics_2}, we get
\begin{equation} \label{eq:fixed_stress_full_estimate1}
\begin{aligned}
    &\left\{\left\lVert \alpha \telta \pfs \right\rVert_{\ltwo}^2-2\alpha\lambda\left(\telta \pfs,\nabla \cdot \telta \ufs\right)+\lambda^2\left\lVert \nabla \cdot \telta \ufs\right\rVert_{\ltwopm}^2 \right\}\\
    &+\left\{(4GD\lambda+\frac{2\gamma D\lambda}{\Delta t})\left\lVert \telta \wfs\right\rVert_{L^2(\mathcal{C})}^2 -2\lambda\left(\telta \pcfs,\telta \wfs\right)_{\mc}+\frac{2\beta_c}{\beta}\left\lVert \telta \pcfs\right\rVert_{L^2(\mc)}^2\right\}\\
    &+\left(\lambda^2+\frac{2\gamma\lambda}{\Delta t}\right)\left\lVert \nabla \cdot \telta \ufs\right\rVert_{\ltwopm}^2+\frac{2\lambda_{\text{min}}\Delta t}{\beta\mu}\left\lVert \nabla\telta \pfs\right\rVert_{\bltwopm}^2+\frac{g_{\text{min}}^3\Delta t}{6\beta\mu}\left\lVert \tnabla\telta  \pcfs\right\rVert_{\bltwoc}^2\\
    &+\frac{g_{\text{min}}}{3\beta\mu}\sum_{*\in\pm}\lVert \telta \pcfs-\telta \pfstr\rVert_{L^2(\mc)}^2\leq \frac{1}{(\beta \lambda)^2}\left\lVert \telta \sigma_{n,k}^h\right\rVert_{\ltwopm}^2+\frac{2}{\beta}\left(\gamma_c\telta \pcpfs-\telta \wpfs,\telta \pcfs\right)_{\mc}.
\end{aligned}
\end{equation}
In view of \eqref{eq:mean_stress}, the quantity inside the first curly brackets of \eqref{eq:fixed_stress_full_estimate1} is equal to $\lVert \telta \sigma_{n,k+1}^{h}\rVert_{\ltwopm}^2$. Moreover, the second term on the right-hand side can be estimated using Young's inequality and picking the coefficient $1/(2\beta_c)$:
\begin{equation}
   \left \lvert \frac{2}{\beta}\left(\gamma_c\telta \pcpfs-\telta \wpfs,\telta \pcfs\right)_{\mc}\right \rvert \leq \frac{1}{\beta\beta_c}\left\lVert \gamma_c \telta \pcpfs-\telta \wpfs\right\rVert_{L^2(\mathcal{C})}^2+\frac{\beta_c}{\beta}\left\lVert \telta \pcfs\right\rVert_{L^2(\mathcal{C})}^2.
\end{equation}
Substituting these into \eqref{eq:fixed_stress_full_estimate1}, we get
\begin{equation} \label{eq:fixed_stress_full_estimate2}
\begin{aligned}
    & \left\lVert \telta \sigma_{n,k+1}^{h} \right\rVert_{\ltwopm}^2+\left\{(4GD\lambda+\frac{2\gamma D\lambda}{\Delta t})\left\lVert \telta \wfs\right\rVert_{L^2(\mathcal{C})}^2 -2\lambda\left(\telta \pcfs,\telta \wfs\right)_{\mc}+\frac{\beta_c}{\beta}\left\lVert \telta \pcfs\right\rVert_{L^2(\mc)}^2\right\}\\
    &+\left(\lambda^2+\frac{2\gamma\lambda}{\Delta t}\right)\left\lVert \nabla \cdot \telta \ufs\right\rVert_{\ltwopm}^2+\frac{2\lambda_{\text{min}}\Delta t}{\beta\mu}\left\lVert \nabla\telta \pfs\right\rVert_{\bltwopm}^2+\frac{g_{\text{min}}^3\Delta t}{6\beta\mu}\left\lVert \tnabla\telta  \pcfs\right\rVert_{\bltwoc}^2\\
    &+\frac{g_{\text{min}}}{3\beta\mu}\sum_{*\in\pm}\lVert \telta \pcfs-\telta \pfstr\rVert_{L^2(\mc)}^2\leq \frac{1}{(\beta \lambda)^2}\left\lVert \telta \sigma_{n,k}^h\right\rVert_{\ltwopm}^2+\frac{1}{\beta \beta_c}\left \lVert \gamma_c \telta \pcpfs - \telta \wpfs\right \rVert_{L^2(\mc)}^2.
\end{aligned}
\end{equation}
It remains to complete the square inside the curly brackets of \eqref{eq:fixed_stress_full_estimate2}. To this end, we recall the definition \eqref{eq:mean_stress_fracture} of the fracture stress, and note that
\begin{equation}
    \left\lVert \telta \xi^{h}_{n,k+1}\right\rVert_{L^2(\mc)}^2=\frac{\gamma_c^2}{\chi^2}\left\lVert \telta \pcfs\right\rVert_{L^2(\mathcal{C})}^2-2\frac{\gamma_c}{\chi^2}\left(\telta \pcfs,\telta \wfs\right)_{\mc}+\frac{1}{\chi^2}\left\lVert \telta \wfs\right\rVert_{L^2(\mc)}^2.
\end{equation}
Comparing with \eqref{eq:fixed_stress_full_estimate2} reveals that the constants $\gamma_c$ and $\chi$ should be chosen such that 
\begin{equation}
    \frac{\beta_c}{\beta}=\frac{\gamma_c^2}{\chi^2}, \ \ \frac{\gamma_c}{\chi^2}=\lambda,
\end{equation}
which gives, using \eqref{eq:splitting_constants}:
\begin{equation}
    \gamma_c=\frac{c_{f}g_{\text{min}}}{\lambda \beta-1}, \ \ \chi=\left(\frac{c_{f}g_{\text{min}}}{\lambda^2\beta-\lambda}\right)^{1/2}.
\end{equation}
We note that $\lambda\beta>1$ and hence, these are well-defined, positive constants. If we now add and subtract $1/\chi^2\lVert \telta \wfs\rVert_{L^2(\mathcal{C})}^2$ to the left-hand side of \eqref{eq:fixed_stress_full_estimate2}, we can complete the square, and the inequality reads:
\begin{equation} \label{eq:final_inequality_proof}
\begin{aligned}
    &\left\lVert \telta \sigma_{n,k+1}^{h}\right\rVert_{\ltwopm}^2+\left\lVert \telta \xi^{h}_{n,k+1}\right\rVert_{L^2(\mathcal{C})}^2+\left(4GD\lambda-\frac{1}{\chi^2}+\frac{2\gamma D\lambda}{\Delta t}\right)\left\lVert \telta \wfs \right\rVert_{L^2(\mathcal{C})}^2\\
    &+\left(\lambda^2+\frac{2\gamma\lambda}{\Delta t}\right)\left\lVert \nabla \cdot \telta \ufs\right\rVert_{\ltwopm}^2+\frac{2\lambda_{\text{min}}\Delta t}{\beta\mu}\left\lVert \nabla\telta \pfs\right\rVert_{\bltwopm}^2+\frac{g_{\text{min}}^3\Delta t}{6\beta\mu}\left\lVert \tnabla\telta \pcfs\right\rVert_{\bltwoc}^2\\
    &+\frac{g_{\text{min}}}{3\beta\mu}\sum_{*\in\pm}\lVert \telta \pcfs-\telta \pfstr\rVert_{L^2(\mc)}^2\leq \frac{1}{(\beta \lambda)^2}\left\lVert \telta \sigma_{n,k}^h\right\rVert_{\ltwopm}^2+\frac{\chi^2}{\beta\beta_c}\left\lVert \telta \xi^{h}_{n,k} \right\rVert_{L^2(\mc)}^2.
\end{aligned}  
\end{equation}
Hence, we end up with two requirements for a contraction:
\begin{equation} \label{eq:contraction_requirements}
    \frac{\chi^2}{\beta\beta_c}<1 \ \ \text{and} \ \  4GD\lambda+\frac{2\gamma D\lambda}{\Delta t}\geq \frac{1}{\chi^2}.
\end{equation}
The first inequality is always satisfied, since by substituting the definitions of $\chi$ and $\beta_c$, it is equivalent to
\begin{equation}
    \frac{1-(\lambda\beta)^2}{\lambda\beta(\lambda\beta-1)}<0,
\end{equation}
which holds since $\lambda\beta>1$. The second inequality of \eqref{eq:contraction_requirements} is equivalent to, after substituting the definition of $\chi$:
\begin{equation} \label{eq:contraction_requirement_2}
    4GD\lambda+\frac{2\gamma D \lambda}{\Delta t}\geq \frac{\lambda^2 c_0}{\alpha^2c_{f}g_{\text{min}}}.
\end{equation}
This inequality can always be satisfied by choosing a sufficiently small time step, and so we assume this. The estimate \eqref{eq:final_inequality} now follows from \eqref{eq:final_inequality_proof} with $C=\max(\frac{1}{(\beta\lambda)^2},\frac{\chi^2}{\beta\beta_c})<1$.
\end{proof}
With the estimate of the preceding Lemma, we can establish the result of Lemma \ref{lemma:fixed_stress_limits}.
\begin{proof}[Proof of Lemma \ref{lemma:fixed_stress_limits}]
Assume a sufficiently small time step $\Delta t$, such that the contraction estimate \eqref{eq:final_inequality} of Lemma \ref{lemma:fixed_stress_contraction} holds. We start by proving the convergence of the fixed-stress iterations $\ppfs$ of the matrix pressure. From the contraction estimate \eqref{eq:final_inequality} it follows that, in particular, $\nabla \telta \ppfs$ converges geometrically to zero, implying that $\nabla \ppfs$ is a Cauchy sequence in $\bltwopm$, and by the Poincaré inequality \eqref{eq:poincare_inequality}, we have that $\ppfs$ is a Cauchy sequence in $\ltwo$. By the completeness of this space, $\ppfs \to p_n^h \in \ltwo$, from which it follows that $\nabla \ppfs \to \nabla p_n^h \in \ltwopm$. Hence, $p_n^h\in \honepm$, and since $\ppfs \in Q^h$ for all $k$, and $Q^h$ is a closed subspace of $\honepm$, we conclude that $p_n^h \in Q^h$.

To show the convergence of the fixed-stress iterations $\pcpfs$ of the fracture pressure, we have from \eqref{eq:final_inequality} that $\tnabla\telta\pcpfs$ converges geometrically to zero, and hence $\tnabla \pcpfs$ is a Cauchy sequence in $\bm{L}^2(\mc)$. Moreover, \eqref{eq:final_inequality} imply that $\xi_{n,k}^h$ and $\wpfs$ are Cauchy sequences in $L^2(\mc)$, and from \eqref{eq:mean_stress_fracture} we then get that $\pcpfs$ is a Cauchy sequence in $L^2(\mc)$. Hence, we conclude that $\pcpfs$ is a Cauchy sequence in $H^1(\mc)$, and by the completeness of this space, we get that $\pcpfs \to \pch \in H^1(\mc)$. Since $\pcpfs\in \Pi^h$ for all $k$, and $\Pi^h$ is a closed subspace of $H^1(\mc)$, we conclude that $\pch \in \Pi^h$.

Finally, we consider the fixed-stress iterations $\upfs$ of the displacement. We first obtain an estimate on $\bm{\varepsilon}(\upfs)$ by applying Young's inequality to \eqref{eq:fixed_stress_estimate_mechanics} and using \eqref{eq:eps_jump_inequality}:
    \begin{equation}
    \begin{aligned}
        &(2G+\frac{\gamma}{\Delta t})\lVert \bm{\varepsilon}(\delta\ufs)\rVert_{\bltwopm}^2+(\lambda+\frac{\gamma}{\Delta t})\lVert \nabla \cdot \delta \ufs\rVert_{\ltwopm}^2 \\
        &\leq \frac{\alpha^2}{2\lambda}\lVert \delta \pfs\rVert_{\ltwo}^2+\frac{\lambda}{2}\lVert \nabla \cdot \delta \ufs\rVert_{\ltwopm}^2+G\lVert \bm{\varepsilon}(\delta\ufs)\rVert_{\bltwopm}^2+\frac{1}{4GD}\lVert\delta \pcfs\rVert_{L^2(\mathcal{C})}^2,
    \end{aligned}
    \end{equation}
    leading to
    \begin{equation}
        (G+\frac{\gamma}{\Delta t})\lVert \bm{\varepsilon}(\delta\ufs)\rVert_{\bltwopm}^2+(\frac{\lambda}{2}+\frac{\gamma}{\Delta t})\lVert \nabla \cdot \delta \ufs\rVert_{\ltwopm}^2\leq\frac{\alpha^2}{2\lambda}\lVert \delta \pfs\rVert_{\ltwo}^2+\frac{1}{4GD}\lVert\delta \pcfs\rVert_{L^2(\mathcal{C})}^2,
    \end{equation}
    and as we have just seen, the right-hand side converges geometrically to zero, implying that $\bm{\varepsilon}(\upfs)$ is a Cauchy sequence in $\bltwopm$. By combining the Poincaré inequality \eqref{eq:poincare_inequality} and the Korn inequality \eqref{eq:korn_inequality}, this implies that $\upfs$ is a Cauchy sequence in $\bhonepm$, and by completeness, $\upfs\to \bm{u}_n^h\in \bhonepm$. Since $\upfs \in \bm{V}^h$ for all $k$, and $\bm{V}^h$ is a closed subspace of $\bhonepm$, we conclude that $\bu_n^h\in \bm{V}^h$.

Turning to problem (Q-D-S), it is straightforward to argue that $p_n^h, \pch$ and $\bu_n^h$ solve the discrete mass balance equations \eqref{eq:discrete_biot} and \eqref{eq:discrete_biot_fracture}, as these equations only involve linear, continuous operators. Hence, by letting $k \to \infty$ in step 1 of the splitting scheme, we immediately recover equations \eqref{eq:discrete_biot} and \eqref{eq:discrete_biot_fracture}. For the variational inequality \eqref{eq:fixed_stress_mechanics} of step 2 of the splitting scheme, the nonlinear contact term that depends on the solution (i.e. the second $J$-functional on the left-hand side) must be considered more carefully. By computing an estimate similar to the one of Lemma \ref{lemma:friction_estimate}, one can show that the contact terms of the variational inequality \eqref{eq:simple_var_ineq} of problem (Q-D-S) is obtained in the limit. For instance, we have for the friction term $J_{\tau}$:
\begin{equation}
\begin{aligned}
&\left \lvert J_{\tau}[\bm{\eta},\delta_t\bu_{n,k}^h]-J_{\tau}[\bm{\eta},\delta_t\bu_n^h] \right \rvert =\\
    &\left \lvert \int_{\mc}\mu_{\text{fr}} c_{\nu}([\bm{\eta}]_{\nu}-g)_+\left\lvert [\delta_t\bu_{n,k}^h]_{\tau}\right\rvert ds-\int_{\mc}\mu_{\text{fr}} c_{\nu}([\bm{\eta}]_{\nu}-g)_+\left\lvert [\delta_t \bu_n^h]_{\tau}\right\rvert ds\right \rvert\\
    &=\left \lvert \int_{\mc}\frac{\mu_{\text{fr}} c_{\nu}}{\Delta t}([\bm{\eta}]_{\nu}-g)_+\left(\vphantom{\rule{0pt}{2.8ex}}\left\lvert [\bm{u}^h_{n,k}-\bm{u}^h_{n-1}]_{\tau}\right\rvert-\left\lvert [\bm{u}^h_n-\bm{u}^h_{n-1}]_{\tau}\right\rvert\right)\right \rvert ds\\
    &\leq \frac{\mu_{\text{fr}} c_{\nu}}{\Delta t}\lVert ([\bm{\eta}]_{\nu}-g)_+ \rVert_{L^2(\mc)}\lVert [\bm{u}^h_{n,k}-\bm{u}^h_n]_{\tau}\rVert_{\bltwoc},
\end{aligned}
\end{equation}
which goes to zero as $k \to \infty$, since $\bm{u}^h_{n,k}\to\bm{u}^h_n$. The same type of estimate can be done on the normal compliance term $\left \lvert J_{\nu}[\bm{\eta},\delta_t\bu_{n,k}^h]-J_{\nu}[\bm{\eta},\delta_t\bu_n^h] \right \rvert$. For the other terms of \eqref{eq:fixed_stress_mechanics}, convergence to the respective terms of \eqref{eq:simple_var_ineq} follows from linearity and continuity. Hence, we conclude that $p_{n}^h, \pch$ and $\bu_{n}^h$ solve problem (Q-D-S).
\end{proof}

\section{Bounding the terms of the error estimate} \label{app:error}

In this appendix, we provide bounds for the terms $A_n^*, B_n^*, C_n^*, D_n^*, E_n^
*, F_n^*, G_n^*, J_n^*$ and $R_n(\dotu_n^I)$ defined in Section \ref{sec:error_estimates}, which are needed to derive the error estimate \eqref{eq:final_estimate}. In the upcoming estimates, we denote generic constants by $C$, while $\varepsilon_i$ denote constants of Young's inequality that are chosen freely (and they must be chosen sufficiently small, as will be shown). With the exception of one case, which will be treated separately, the constants involved are independent of the discretization parameters.

Starting with the $A_n^*$-term, we use the continuity of the bilinear form \eqref{eq:continuity_bilinear} and the triangle inequality to obtain
\begin{equation} \label{eq:error_comp_1}
    \begin{aligned}
        \Delta t\sum_{n=1}^N \lvert A_n^* \rvert
        \leq m_A\Delta t\sum_{n=1}^N &\left(\vphantom{\rule{0pt}{3ex}}\lVert \dot{a}_n^{\bu}\rVert_{\bhonepm}\left(\lVert \inu \rVert_{\bhonepm}+\lVert \enu \rVert_{\bhonepm}\right)\right.\\
        &+\lVert \enu \rVert_{\bhonepm}\left( \lVert \dotu_n-\delta_t\bu_n\rVert_{\bhonepm}+\lVert \delta_t \inu \rVert_{\bhonepm}\right)\\
        &\left.+\lVert \inu \rVert_{\bhonepm}\left( \lVert \dotu_n-\delta_t\bu_n\rVert_{\bhonepm}+\lVert \delta_t \inu \rVert_{\bhonepm}+\lVert \delta_t \enu \rVert_{\bhonepm}\right)\vphantom{\rule{0pt}{3ex}}\right).\\
    \end{aligned}
\end{equation}
Using Young's inequality on each term on the right-hand side of \eqref{eq:error_comp_1} yields, after merging the various constants:
\begin{equation} \label{eq:error_comp_2}
\begin{aligned}
    \Delta t \sum_{n=1}^N \lvert A_n^* \rvert &\leq \varepsilon_1\sum_{n=1}^N \Delta t \lVert \delta_t \enu \rVert_{\bhonepm}^2+\varepsilon_2\Delta t \sum_{n=1}^N \lVert \enu \rVert_{\bhonepm}^2\\
    &+C\Delta t \sum_{n=1}^N \left(\lVert \dotu_n-\delta_t \bu_n \rVert_{\bhonepm}^2+\lVert \delta_t \inu \rVert_{\bhonepm}^2+\lVert \inu \rVert_{\bhonepm}^2+\lVert \dot{a}_n^{\bu}\rVert_{\bhonepm}^2\right).
\end{aligned}
\end{equation}
The first term inside the third sum on the right-hand side of \eqref{eq:error_comp_2} can be bounded, using Taylor's theorem with integral remainder, as follows:
\begin{equation} \label{eq:taylor_estimate_1}
  \Delta t \sum_{n=1}^N\lVert \dotu_n-\frun\rVert_{\bhonepm}^2 \leq C(\Delta t)^2 \lVert \ddot{\bu}\rVert_{L^2(0,T;\bhonepm)}^2,  
\end{equation}
while the second term inside the third sum on the right-hand side of \eqref{eq:error_comp_2} can be bounded, using the Fundamental Theorem of Calculus, as follows:
\begin{equation} \label{eq:taylor_estimate_2}
    \Delta t\sum_{n=1}^N \lVert \franu \rVert_{\bhonepm}^2=\sum_{n=1}^N \frac{1}{\Delta t}\lVert \inu-a_{n-1}^{\bu} \rVert_{\bhonepm}^2 \leq \lVert \dotu-\dotu^I\rVert_{L^2(0,T;\bhonepm)}^2,
\end{equation}
Hence, the final estimate for the $A_n^*$-term becomes:
\begin{equation} \label{eq:estimate_A}
    \begin{aligned}
            \Delta t \sum_{n=1}^N \lvert A_n^* \rvert &\leq \varepsilon_1\sum_{n=1}^N \Delta t \lVert \delta_t \enu \rVert_{\bhonepm}^2+\varepsilon_2\Delta t \sum_{n=1}^N \lVert \enu \rVert_{\bhonepm}^2+C\Bigg(\lVert \dotu-\dotu^I \rVert_{L^2(0,T;\bhonepm)}^2\\
            &+(\Delta t)^2\lVert \ddot{\bu} \rVert_{L^2(0,T;\bhonepm)}^2+\Delta t\sum_{n=1}^N \left(\lVert \inu \rVert_{\bhonepm}^2+\lVert \dot{a}_n^{\bu}\rVert_{\bhonepm}^2\right)\Bigg).
    \end{aligned}
\end{equation}
The terms $B_n^*-G_n^*$ are estimated using mostly the same techniques as for the $A_n^*$-term, with only some minor differences. Therefore, instead of showing the full derivation of each of these estimates, we will merely state the final result, and point out the minor differences. 
The $B_n^*$-estimate is derived by the exact same techniques as the $A_n^*$-estimate, and the result is:
\begin{equation} \label{eq:estimate_B}
    \begin{aligned}
        \Delta t \sum_{n=1}^N \lvert B_n^* \rvert \leq \varepsilon_3\Delta t \sum_{n=1}^N \lVert \delta_t \enu\rVert_{\bhonepm}^2+C\left(\vphantom{\sum_1^N}\lVert \dotu-\dotu^I\rVert^2_{L^2(0,T;\bhonepm}\right.&+(\Delta t)^2\lVert \ddot{\bu} \rVert_{L^2(0,T;\bhonepm)}^2\\
        &\left.+\Delta t\sum_{n=1}^N\lVert \dot{a}_n^{\bu}\rVert_{\bhonepm}^2\right),
    \end{aligned}
\end{equation}
For the $C_n^*$-term, in addition to the previous techniques, we use the fact that $\lVert \nabla \cdot \bm{v}\rVert_{\ltwopm} \leq C\lVert \bm{v} \rVert_{\bhonepm}$, for any $\bm{v} \in \bhonepm$, for some constant $C$ independent of discretization parameters. The result is:
\begin{equation} \label{eq:estimate_C}
    \begin{aligned}
        \Delta t \sum_{n=1}^N \lvert C_n^* \rvert &\leq \varepsilon_4\Delta t\sum_{n=1}^N \lVert\delta_t \enu \rVert_{\bhonepm}^2+\varepsilon_5\Delta t\sum_{n=1}^N \lVert \enp \rVert_{\ltwo}^2+C\Bigg(\lVert \dotu-\dotu^I\rVert^2_{L^2(0,T;\bhonepm)}\\
        &+(\Delta t)^2\lVert \ddot{\bu} \rVert_{L^2(0,T;\bhonepm)}^2+\Delta t\sum_{n=1}^N \left( \lVert \inp \rVert_{\ltwo}^2+\lVert \dot{a}_n^{\bu}\rVert_{\bhonepm}^2\right)\Bigg),
    \end{aligned}
\end{equation}
For the $D_n^*$-term, we use the trace inequality \eqref{eq:trace_inequality} to bound the jump terms, in addition to the previous techniques. The result is:
\begin{equation} \label{eq:estimate_D}
\begin{aligned}
    \Delta t\sum_{n=1}^N \lvert D_n^*\rvert &\leq \varepsilon_6\Delta t \sum_{n=1}^N \lVert \frenu \rVert_{\bhonepm}^2+\varepsilon_7\Delta t \sum_{n=1}^N \lVert \enpc \rVert_{L^2(\mc)}^2+C\Bigg( \lVert \dotu-\dotu^I\rVert^2_{L^2(0,T;\bhonepm)}\\ 
    &+(\Delta t)^2\lVert \ddot{\bu} \rVert_{L^2(0,T;\bhonepm)}^2+\Delta t \sum_{n=1}^N \big(\lVert \inpc \rVert_{L^2(\mc)}^2+\lVert \dot{a}_n^{\bu}\rVert_{\bhonepm}^2\big)\Bigg)
\end{aligned}
\end{equation}
For the $E_n^*-$ and $F_n^*$-terms, we use formulas analogous to \eqref{eq:taylor_estimate_1} and \eqref{eq:taylor_estimate_2} to bound terms of the form $\Delta t \sum_{n=1}^N\lVert \delta_t \inp \rVert_{\ltwo}^2$, $\Delta t \sum_{n=1}^N \lVert \delta_t \inpc \rVert_{L^2(\mc)}^2$, $\Delta t \sum_{n=1}^N\lVert \dotp_n-\delta_t p_n \rVert_{\ltwo}^2$ and $\Delta t \sum_{n=1}^N \lVert \dotpc_n-\delta_t \pcn \rVert_{L^2(\mc)}^2$. Additionally, we use the boundedness estimates \eqref{eq:spd_perm} and \eqref{eq:g_bound} for the permeability $\bk$ and fracture gap $g$, respectively, and absorb the corresponding constants $\lambda_{\text{max}}, g_{\text{max}}$ into the other constants. The resulting estimates become:
\begin{equation} \label{eq:estimate_F}
    \begin{aligned}
        &\Delta t \sum_{n=1}^N \lvert E_n^* \rvert \leq \varepsilon_8\Delta t\sum_{n=1}^N\lVert \enp \rVert_{\ltwo}^2+\varepsilon_9\Delta t \sum_{n=1}^N \lVert \nabla \enp \rVert_{\bltwopm}^2+C\Bigg(\lVert \dotp-\dotp^I \rVert_{L^2(0,T;\ltwo)}^2\\
        &+\lVert \dotu-\dotu^I \rVert_{L^2(0,T;\bhonepm)}^2+(\Delta t)^2\left(\lVert \ddot{p}\rVert_{L^2(0,T;\ltwo)}^2+\lVert \ddot{\bu} \rVert_{L^2(0,T;\bhonepm)}^2\right)+\Delta t \sum_{n=1}^N \lVert \nabla \inp \rVert_{\bltwopm}^2\Bigg).
    \end{aligned}
\end{equation}
\begin{equation}
    \begin{aligned}
        &\Delta t \sum_{n=1}^N \lvert F_n^* \rvert \leq \varepsilon_{10}\Delta t \sum_{n=1}^N \lVert \enpc \rVert_{L^2(\mc)}^2+\varepsilon_{11}\Delta t \sum_{n=1}^N \lVert \tnabla\enpc\rVert_{L^2(\mc)}^2+C\Bigg(\lVert \dotpc-\dotpc^I \rVert_{L^2(0,T;L^2(\mc))}^2\\
        &+\lVert \dotu-\dotu^I \rVert_{L^2(0,T;\bhonepm)}^2+(\Delta t)^2\Big(\lVert \ddotpc\rVert_{L^2(0,T;L^2(\mc))}^2+\lVert \ddot{\bu}\rVert_{L^2(0,T;\bhonepm)}^2\Big)+\Delta t \sum_{n=1}^N\lVert \tnabla \inpc\rVert_{L^2(\mc)}^2\Bigg).
    \end{aligned}
\end{equation}
For the $G_n^*$-estimate, we use the trace inequality \eqref{eq:trace_inequality} to obtain the bound $\sum_{*\in\pm}\lVert \inptr \rVert_{L^2(\mc)}^2 \leq C \lVert \inp \rVert_{\honepm}^2$. We also use once again the boundedness of the fracture gap, and absorb $g_{\text{max}}$ into the other constants, to get the estimate:
\begin{equation} \label{eq:estimate_E}
    \begin{aligned}
        \Delta t\sum_{n=1}^N\lvert G_n^*\rvert \leq \varepsilon_{12}\sum_{n=1}^N\Delta t\sum_{*\in\pm}\lVert \enpc-\enptr\rVert_{L^2(\mc)}^2+C\Bigg(\Delta t \sum_{n=1}^N\big(\lVert \inpc\rVert_{L^2(\mc)}^2+\lVert \inp\rVert_{\honepm}^2\big)\Bigg),
    \end{aligned}
\end{equation}
Next, we consider the $J_n^*$-term. It is bounded by Lemma \ref{lemma:friction_estimate}, and by writing $\dotu_n^I-\frunh=\dotu_n^I-\dotu_n+\dotu_n-\frun+\frun-\frunh$: 
\begin{equation} \label{eq:j_estimate}
\begin{aligned}
    \lvert J_n^*\rvert &= \left\lvert J[\bu_n,\frunh]-J[\bu_n,\dotu_n^I]+J[\bu_n^h,\dotu_n^I]-J[\bu_n^h,\frunh] \right\rvert\\
    &\leq C_J \lVert \bu_n-\bu_n^h \rVert_{\bhonepm}\lVert \dotu_n^I-\frunh \rVert_{\bhonepm}\\
    &\leq C_J\lVert \inu+\enu \rVert_{\bhonepm}(\lVert \dotu_n^I-\dotu_n\rVert_{\bhonepm}+\lVert \dotu_n-\frun\rVert_{\bhonepm}+\lVert \frun-\frunh\rVert_{\bhonepm})\\
    &= C_J\lVert \inu+\enu \rVert_{\bhonepm}(\lVert \dot{a}_n^{\bu}\rVert_{\bhonepm}+\lVert \dotu_n-\frun\rVert_{\bhonepm}+\lVert \franu+\frenu\rVert_{\bhonepm}).
\end{aligned}
\end{equation} 
From here we perform the same kinds of computations as for the previous terms to get
\begin{equation} \label{eq:estimate_J}
    \begin{aligned}
       &\Delta t \sum_{n=1}^N \lvert J_n^* \rvert \leq (\varepsilon_{13}+\varepsilon_{14}) \sum_{n=1}^N\Delta t \lVert \delta_t \enu \rVert_{\bhonepm}^2+((4\varepsilon_{13})^{-1}+\varepsilon_{15})\sum_{n=1}^N \Delta t \lVert \enu \rVert_{\bhonepm}^2\\
       &+C\left(\vphantom{\sum_1^N}\lVert \dotu-\dotu^I \rVert_{L^2(0,T;\bhonepm)}^2+(\Delta t)^2\lVert \ddot{\bu} \rVert_{L^2(0,T;\bhonepm)}^2\right.\left.+\Delta t\sum_{n=1}^N\left(\lVert \inu \rVert_{\bhonepm}^2+\lVert \dot{a}_n^{\bu}\rVert_{\bhonepm}^2\right)\right).
    \end{aligned}
\end{equation}
Here we get the same caveat as in the stability estimate, namely that we obtain the term \\ $(4\varepsilon_{13})^{-1}\sum_{n=1}^N\Delta t\lVert \enu \rVert_{\bhonepm}^2$ whose size can only be controlled by imposing a restriction on the time step size. This leads to the same restriction as in the stability estimate, $\Delta t<2\alpha_A\alpha_B$, required to employ Grönwall's inequality at the end of the proof.

Finally, we consider the residual term $R_n$. By simple usage of Cauchy-Schwarz, continuity of the bilinear forms, and an estimate on the $J$-terms similar to Lemma \ref{lemma:friction_estimate}, we can bound the term by
\begin{equation}
    \begin{aligned}
        \lvert R_n(\dotu_n^I)\rvert &\leq \left\lvert A[\bu_n,\dot{a}_n^{\bu}]\right\rvert+\left\lvert B[\dotu_n,\dot{a}_n^{\bu}]\right\rvert+\alpha\left\lvert\left(p_n,\nabla \cdot \dot{a}_n^{\bu}\right)\right\rvert+\lvert(\pcn,[\dot{a}_n^{\bu}]_{\nu})_{\mc}\rvert\\
        &\ \ \ \ \ \ \ \ \ \ \ \ \ \ \ \ \ \ \ \ \ \ \ \ \ \ \ \ \ \ \ \ \ \ \ +\left\lvert J[\bu_n,\dotu_n^I]-J[\bu_n,\dotu_n]\right\rvert+\left\lvert\left(\bm{f}_n,\dot{a}_n^{\bu}\right)\right\rvert\\
        &\leq C_n \lVert \dot{a}_n^{\bu}\rVert_{\bhonepm},
    \end{aligned}
\end{equation}
where $C_n$ depends on the (continuous) solution and data at time step $t_n$, and hence, it is dependent on the temporal discretization. This is in contrast to the previous constants, which have all been independent of discretization parameters. However, our assumptions on the temporal regularity of the solution and the data are sufficient to ensure that $C_n$ is bounded on $[0,T]$, and we can define $C$ to be any constant such that $\max_{1\leq n \leq N}C_n<C$ for all $N$, and hence obtain the estimate:
\begin{equation}
    \Delta t \sum_{n=1}^N \lvert R_n(\dotu_n^I)\rvert \leq C \Delta t \sum_{n=1}^N \lVert \dot{a}_n^{\bu}\rVert_{\bhonepm},
\end{equation}
Note that this term does not appear in squared form, which will result in a reduced spatial convergence order.

The $\varepsilon_i$ constants are now chosen such that the corresponding terms can be absorbed into the left-hand side of \eqref{eq:big_estimate_summed}. We accordingly choose $\varepsilon_1+\varepsilon_3+\varepsilon_4+\varepsilon_6+\varepsilon_{13}+\varepsilon_{14}<\alpha_B$, $\varepsilon_{12}<g_{\text{min}}/(6\mu)$, $\varepsilon_{9}<\lambda_{\text{min}}/\mu$ and $\varepsilon_{11}<g_{\text{min}}^3/(12\mu)$. Finally, assuming the aforementioned restriction on the time step size, we employ Grönwall's inequality on the terms $((4\varepsilon_{13})^{-1}+\varepsilon_2+\varepsilon_{15})\Delta t \sum_{n=1}^N \lVert \enu \rVert_{\bhonepm}^2$, $(\varepsilon_5+\varepsilon_8)\Delta t\sum_{n=1}^N \lVert \enp \rVert_{\ltwo}^2$ and $(\varepsilon_7+\varepsilon_{10})\Delta t \sum_{n=1}^N\lVert \enpc \rVert_{L^2(\mc)}^2$. After combining the various constants, we arrive at the estimate \eqref{eq:final_estimate}.

\section*{Data Availability}

The simulation results of this paper can be reproduced using a Docker image available at Zenodo \cite{zenodo_fractured_biot}.

\printbibliography

\end{document}